\documentstyle[newlfont,11pt]{amsart}

\newcommand{\nc}{\newcommand}

\nc{\one}{\mbox{\bf 1}}
\nc{\invtensor}{\underset{\leftarrow}{\otimes}}
\nc{\rlarrows}{\begin{picture}(1,0.4)
                \put(0,-0.1){\makebox(1,0.2){$\leftarrow$}}
                \put(0,0.1){\makebox(1,0.2){$\ra$}}
              \end{picture}}
\nc{\rra}{\begin{picture}(1,0.4)
                \put(0,-0.1){\makebox(1,0.2){$\lra$}}
                \put(0,0.1){\makebox(1,0.2){$\lra$}}
              \end{picture}}
\nc{\Left}{\bold L}
\nc{\Right}{\bold R}
\nc{\gr}{\operatorname{gr}}
\nc{\Ker}{\operatorname{Ker}}
\nc{\Ho}{\operatorname{Ho}}
\nc{\alt}{\operatorname{alt}}
\nc{\Sym}{\operatorname{Sym}}
\nc{\sym}{\operatorname{sym}}
\nc{\id}{\operatorname{id}}
\nc{\Der}{\operatorname{Der}}
\nc{\im}{\operatorname{Im}}
\nc{\Col}{\operatorname{Col}}
\nc{\ter}{\operatorname{ter}}
\nc{\intl}{\operatorname{int}}
\nc{\out}{\operatorname{out}}
\nc{\irr}{\operatorname{irr}}
\nc{\Tor}{\operatorname{Tor}}
\nc{\WE}{\frak W}
\nc{\FIB}{\frak F}
\nc{\COF}{\frak C}
\nc{\Sh}{\operatorname{Sh}}
\nc{\Hom}{\operatorname{Hom}}
\nc{\CHom}{{\cal H}om}
\nc{\End}{\operatorname{End}}
\nc{\holim}{\operatorname{holim}}
\nc{\dirlim}{\underset{\rightarrow}{\lim}\,}
\nc{\invlim}{\underset{\leftarrow}{\lim}\,}
\nc{\CB}{\operatorname{\bf CB}}
\nc{\com}{\operatorname{co}}
\nc{\Tot}{\operatorname{Tot}}
\nc{\Th}{\operatorname{Th}}
\nc{\Cech}{\check{C}}
\nc{\Spec}{\operatorname{Spec}}
\nc{\MC}{\operatorname{MC}}
\nc{\U}{\operatorname{U}}
\nc{\Diff}{{\cal D}\mbox{\em iff}}
\nc{\Mor} {\operatorname{Mor}}
\nc{\Ob}{\operatorname{Ob}}
\nc{\Cone}{\operatorname{cone}}
\nc{\Coder}{\operatorname{Coder}}
\nc{\cone}{\operatorname{cone}}
\nc{\Liee}{\operatorname{Lie}}
\nc{\dbar}{\bar{\partial}}
\nc{\tr}{\operatorname{tr}}
\nc{\is}{\operatorname{iso}}

\renewcommand{\mod}{\mathtt {Mod}}          
\nc{\Mod}{\mathtt {Mod}}
\nc{\MOR}{\mathtt {MOR}}
\nc{\Alg}{\mathtt {Alg}}
\nc{\dgl}{\mathtt {DGL}}
\nc{\dga}{\mathtt {DGA}}
\nc{\cdga}{\mathtt {DGC}}
\nc{\hol}{\mathtt {Holie}}
\nc{\hoass}{\mathtt {Hoass}}
\nc{\hocom}{\mathtt {Hocom}}
\nc{\hoa}{\mathtt {Hoalg}}
\nc{\hoalg}{\mathtt {Hoalg}}
\nc{\Ens}{\mathtt {Ens}}
\nc{\Cat}{\mathtt {Cat}}

\nc{\Ass}{\mathtt {Ass}}
\nc{\Com}{\mathtt {Com}}
\nc{\Lie}{\mathtt {Lie}}
\nc{\op}{\mathtt {Op}}
\nc{\Op}{\mathtt {Op}}
\nc{\asop}{\mathtt {Asop}}

\nc{\pa}{\partial}

\nc{\CA}{\cal A}
\nc{\CC}{\cal C}
\nc{\CD}{\cal D}
\nc{\CF}{\cal F}
\nc{\CI}{\cal I}
\nc{\CJ}{\cal J}
\nc{\CK}{\cal K}
\nc{\CL}{\cal L}
\nc{\CO}{\cal O}
\nc{\CU}{\cal U}
\nc{\CS}{\cal S}
\nc{\CT}{\cal T}

\nc{\fg}{\frak g}
\nc{\fh}{\frak h}
\nc{\fs}{\frak s}
\nc{\ft}{\frak t}

\nc{\nen}{\newenvironment}
\nc{\ol}{\overline}
\nc{\ul}{\underline}
\nc{\ra}{\rightarrow}
\nc{\lra}{\longrightarrow}
\nc{\lla}{\longleftarrow}
\nc{\Lra}{\Longrightarrow}
\nc{\Lla}{\Longleftarrow}
\nc{\Llra}{\Longleftrightarrow}
\nc{\hra}{\hookrightarrow}
\nc{\iso}{\overset{~}{\lra}}

\nc{\Thm}[1]{Theorem~\ref{#1}}
\nc{\Prop}[1]{Proposition~\ref{#1}}
\nc{\Lem}[1]{Lemma~\ref{#1}}
\nc{\Cor}[1]{Corollary~\ref{#1}}
\nc{\Conj}[1]{Conjecture~\ref{#1}}
\nc{\Claim}[1]{Claim~\ref{#1}}
\nc{\Defn}[1]{Definition~\ref{#1}}
\nc{\Exa}[1]{Example~\ref{#1}}
\nc{\Rem}[1]{Remark~\ref{#1}}
\nc{\Note}[1]{Note~\ref{#1}}
\nc{\Sect}[1]{Section~\ref{#1}}

\nen{thm}[1]{\label{#1}{\bf Theorem.\ } \em}{}
\nen{prop}[1]{\label{#1}{\bf Proposition.\ } \em}{}
\nen{lem}[1]{\label{#1}{\bf Lemma.\ } \em}{}
\nen{cor}[1]{\label{#1}{\bf Corollary.\ } \em}{}
\nen{conj}[1]{\label{#1}{\bf Conjecture.\ } \em}{}
\nen{claim}[1]{\label{#1}{\bf Claim.\ } \em}{}

\nen{defn}[1]{\label{#1}{\bf Definition.\ } }{}
\nen{exa}[1]{\label{#1}{\bf Example.\ } }{}

\nen{rem}[1]{\label{#1}{\em Remark.\ } }{}
\nen{note}[1]{\label{#1}{\em Note.\ } }{}
\nen{exer}[1]{\label{#1}{\em Exercise.\ } }{}

\begin{document}

\title[]{Homological algebra of homotopy algebras}
\author{Vladimir Hinich}
\address{Dept. of Mathematics and Computer Science, University of Haifa,
Mount Carmel, Haifa 31905 Israel}



\maketitle

\tableofcontents
\section{Introduction}

\subsection{}
Let $k$ be a ring and let $C(k)$ be the category of (unbounded) complexes
of $k$-modules. As we know from Spaltenstein~\cite{sp} ( see 
also~\cite{a} and~\cite{belu})  one can do a homological
algebra in $C(k)$ using the appropriate notions for $\CK$-projective and
$\CK$-injective complexes. 

The present paper started from the observation that this homological algebra
in $C(k)$ (or, more generally, in the category of dg modules over an
associative dg algebra) can be described using Quillen's language of closed 
model categories (see~\cite{q1,q2}). 
For, if we take quasi-isomorphisms in $C(k)$ to be weak equivalences and 
componentwise surjective maps of complexes to be 
fibrations, then a closed model category structure on $C(k)$ is defined, and 
cofibrant objects in it are precisely the $\CK$-projectives of~\cite{sp}. 

The possibility of working with unbounded complexes is very important if
we wish to work with ``weak algebras'' --- the ones satisfying the 
standard identities (associativity, commutativity, or Jacobi identity,
for example) up to some higher homotopies. An appropriate language to
describe these objects is that of operads (see~\cite{km} and 
references therein) and one feels extremely incomfortable when 
restricted to, say, non-negatively graded world (for instance
non-negatively graded commutative dg algebras do not admit semi-free
resolutions; enveloping algebras of operad algebras are very often 
infinite in both directions). 

In this paper we use Quillen's machinery of closed model categories
to describe homological algebra connected to operads and operad
algebras.

In Sections~\ref{cmc}---~\ref{section-catop} we define the necessary
structures and prove some standard comparison results. In 
Section~\ref{section-cot} we study the notion of
cotangent complex of a morphism of dg operad algebras. In the last
Section~\ref{section-tan} we define a canonical
structure of homotopy Lie algebra on the tangent complex. The latter 
is the main (concrete) result of the paper.  

Let us describe in a more detail the contents of the paper.

\subsection{}
Homological algebra of operad algebras has three different levels.

On the lowest level we have the category $\mod(\CO,A)$ of modules over a 
fixed algebra $A$ over an operad $\CO$. This is the category of dg modules
over the enveloping algebra $U(\CO,A)$ which is an associative dg algebra.
As we mentioned above, this category admits a closed model category (CMC)
structure --- see~\ref{dha1-models}; the corresponding homotopy category
is the derived category of $U(\CO,A)$-modules and it is denoted by 
$DU(\CO,A)$. 

Since operad algebra $(\CO,A)$ in $C(k)$ is not just a firm collection 
of operations but is merely a model for the idea of "algebra up to 
homotopy", we have to understand what happens to $DU(\CO,A)$ when one
substitutes $(\CO,A)$ with a quasi-isomorphic algebra $(\CO',A')$.

On the next level we have the category $\Alg(\CO)$ of algebras over a 
fixed operad $\CO$. This category also admits a CMC structure, provided
some extra hypotheses on $\CO$ ($\Sigma$-splitness, see~\ref{ssplit}) are 
fulfilled. These extra hypotheses correspond more or less to the cases
where one is able to use free algebra-resolutions instead of simplicial
resolutions, in order to define algebra cohomology (see~\cite{q3}). Thus,
the operad $\Ass_k$ responsible for associative $k$-algebras is $\Sigma$-split
for any $k$; any operad over $k$ is $\Sigma$-split when $k\supseteq\Bbb Q$.

Finally, on the highest level we have the category $\op(k)$ of (dg)
operads over $k$. Quasi-isomorphic operads here correspond, roughly
speaking, to different collections of higher homotopies used in an
algebra $A$ in order to make it ``homotopy algebra''.
The category $\op(k)$ also admits a CMC structure. 

What is the connection between the different model structures?

First of all, if one has a quasi-isomorphism $\alpha:\CO\to\CO'$
of $\Sigma$-split operads compatible with $\Sigma$-splitting
(this condition is fulfilled, e.g., when $k\supseteq\Bbb{Q}$) then
the homotopy categories $\hoa(\CO)$ and $\hoa(\CO')$ are naturally
equivalent --- see~\Thm{equivalence}. This result implies, for instance, the
representability of strong homotopy algebras (Lie or not) by
strict algebras in characteristic zero.

A similar equivalence on the lower level takes place only for
associative algebras: \Thm{asso-comparison} claims that a 
quasi-isomorphism $f:A\to B$ of associative dg algebras induces
an equivalence of the derived categories $D(A)$ and $D(B)$.
For algebras over an arbitrary operad $\CO$ one has such a comparison
result only when $A$ and $B$ are cofibrant algebras (see~\Cor{U(w-c)}),
or when the operad $\CO$ is cofibrant and $A, B$ are flat as 
$k$-complexes. This suggests a definition of derived category $D(\CO,A)$
which can be "calculated" either as $D(\CO,P)$ where $P$ is a cofibrant 
$\CO$-algebra quasi-isomorphic to $A$, or as $D(\tilde{\CO},A)$
where $\tilde{\CO}$ is a cofibrant resolution of $\CO$ and $A$ is flat. 
This is done in~\ref{virtual} and in~\ref{more}. The category
$D(\CO,A)$ is called the derived category of virtual $A$-modules and it
depends functorially on $(\CO,A)$.

\subsection{} For any morphism $f:A\to B$ of $\CO$-algebras one defines
in a standard way the functor $\Der$:
$$ \Der_{B/A}:\Mod(\CO,B)\to C(k).$$
The functor is representable by the module $\Omega_{B/A}\in\Mod(\CO,B)$.
This is the module of differentials. If $\CO$ is $\Sigma$-split 
so that $\Alg(\CO)$ admits a CMC structure, one defines the relative
cotangent complex $L_{B/A}\in D(\CO,B)$ as the module of differentials
of a corresponding cofibrant resolution. This defines cohomology
of $\CO$-algebra $A$ as the functor
$$ M\in D(\CO,A)\mapsto H(A,M)=\Right\Hom(L_A,M)\in C(k).$$

The most interesting cohomology is the one with coefficients in $M=A$.
No doubt, the complex $\Der^{\CO}(A,A)$ admits a dg Lie algebra structure.
The main result of Section~\ref{section-tan}, \Thm{T-is-a-functor}, claims that
the tangent complex $T_A:=H(A,A)=\Right\Hom(L_A,A)$ admits a canonical
structure of Homotopy Lie algebra. This means that $T_A$ is defined uniquely
up to a unique isomorphism as an object of the category $\hol(k)$. 

\subsection{} Let us indicate some relevant references.

Spaltenstein~\cite{sp}, Avramov-Foxby-Halperin~\cite{a} developed
homological algebra for unbounded complexes. 

Operads and operad algebras
were invented by J.P.~May in early 70-ies in a topological context; 
dg operads appeared explicitly in~\cite{hlha} and became popular in 90-ies 
mainly because of their connection to quantum field theory.

M.~Markl in~\cite{mar} studied ``minimal models'' for operads --- similarly
to Sullivan's minimal models for commutative dg algebras over $\Bbb Q$.
In our terms, these are cofibrant operads weakly equivalent to a given one.

In~\cite{ss} M..~Schlessinger and J.~Stasheff propose to define the tangent
complex of a commutative algebra $A$ as $\Der(\CA)$ where $\CA$
is a ``model'' i.e. a commutative dg algebra quasi-isomorphic to $A$ and 
free as a graded commutative algebra. This complex has an obvious Lie
algebra structure which is proven to coincide sometimes (for a standard
choice of $\CA$) with the one defined by the Harrison complex of $A$. 

It is clear ``morally'' that the homotopy type of the Lie algebra $\Der(\CA)$
should not depend on the choice of $\CA$.
Our main result of Section~\ref{section-tan} says (in a more general setting)
that this is really so.  

\subsection{Notations} For a ring $k$ we denote by $C(k)$ the category of
complexes of $k$-modules. If $X,Y\in C(k)$ we denote by $\CHom_k(X,Y)$
the complex of maps form $X$ to $Y$ (not necessarily commuting with
the differentials).

$\Bbb N$ is the set of non-negative integers; 
$\Ens$ is the category of sets, $\Cat$ is the $2$-category of small categories.
The rest of the notations is given in the main text.

\section{Closed model categories}
\label{cmc}

The main result of this Section --- \Thm{cmc-structure} --- provides 
a category $\CC$ endowed with a couple of adjoint functors
$$ \# :\CC\rlarrows C(k): F$$
($F$ is left adjoint to $\#$) where $C(k)$ is the category of (unbounded)
complexes of modules over a ring $k$  and satisfying properties (H0), 
(H1) of~\ref{first}, with a structure of closed model category (CMC)
in sense of Quillen \cite{q1},\cite{q2}, see also ~\ref{recall-cmc}. 
This allows one to define a CMC structure on the
category of (dg) operad algebras (\Sect{section-opalg}), on the
category of modules over an associative dg algebra  (\Sect{dha-1})
and, more generally, on the category of modules over an operad algebra
(\Sect{moa}). The CMC structure on the category of operads 
(\Sect{section-catop}) is obtained in almost the same way.

\subsection{Definition}
\label{recall-cmc}
Recall (cf.~\cite{q1},~\cite{q2}) that a closed model category
(CMC) structure on a category $\CC$ is given by three collections of
morphisms ---
{\em weak equivalences} ($\WE$),  {\em fibrations} ($\FIB$),  
{\em cofibrations} ($\COF$) in $\Mor(\CC)$ such that the following
axioms are fulfilled:

(CM 1) $\CC$ is closed under finite limits and colimits.

(CM 2) Let $f,g\in\Mor(\CC)$ such that $gf$ is defined. If any two
of $f,g,gf$ are in $\WE$ than so is the third one.

(CM 3) Suppose that $f$ is a retract of $g$ i.e. that there exists
a commutative diagram
\begin{center}
\begin{picture}(5,3)
   \put(0,0){\makebox(1,1){$\bullet$}}
   \put(2,0){\makebox(1,1){$\bullet$}}
   \put(4,0){\makebox(1,1){$\bullet$}}
   \put(0,2){\makebox(1,1){$\bullet$}}
   \put(2,2){\makebox(1,1){$\bullet$}}
   \put(4,2){\makebox(1,1){$\bullet$}}

   \put(0.5,2){\vector(0,-1){1}}
   \put(2.5,2){\vector(0,-1){1}}
   \put(4.5,2){\vector(0,-1){1}}
   \put(1,0.5){\vector(1,0){1}}
   \put(1,2.5){\vector(1,0){1}}
   \put(3,0.5){\vector(1,0){1}}
   \put(3,2.5){\vector(1,0){1}}

   \put(0,1){\makebox(0.5,1){$\scriptstyle f$}}
   \put(2,1){\makebox(0.5,1){$\scriptstyle g$}}
   \put(4,1){\makebox(0.5,1){$\scriptstyle f$}}

\end{picture}
\end{center}

in which the compositions of the horisontal maps are identities.
Then if $g$ belongs to $\WE$ (resp., $\FIB$ or $\COF$) then so does 
$f$.

(CM 4) Let
\begin{center}
\begin{picture}(3,3)
   \put(0,0){\makebox(1,1){$B$}}
   \put(2,0){\makebox(1,1){$Y$}}
   \put(0,2){\makebox(1,1){$A$}}
   \put(2,2){\makebox(1,1){$X$}}
   \put(0.5,2){\vector(0,-1){1}}
   \put(2.5,2){\vector(0,-1){1}}
   \put(1,2.5){\vector(1,0){1}}
   \put(1,0.5){\vector(1,0){1}}

   \multiput(0.8,0.8)(0.1,0.1){14}{\makebox(0.1,0.1){$\cdot$}}
   \put(2.2,2.2){\vector(1,1){0}}

   \put(0,1){\makebox(0.5,1){$\scriptstyle i$}}
   \put(2.5,1){\makebox(0.5,1){$\scriptstyle p$}}
   \put(1,1){\makebox(1.15,0.65){$\scriptstyle\alpha$}}

\end{picture}
\end{center}

be a commutative diagram with $i\in\COF,\ p\in\FIB$. Then a dotted
arrow $\alpha$ making the diagram commutative, exists if
either

(i) $i\in\WE$  

or

(ii) $p\in\WE$.

(CM 5) Any map $f:X\ra Y$ can be decomposed in the following ways:

(i) $f=pi$, $p\in\FIB,\ i\in\WE\cap\COF$;

(ii) $f=qj$, $q\in\WE\cap\FIB,\ j\in\COF$.

The morphisms in $\WE\cap\FIB$ are called {\em acyclic fibrations};
the morphisms in $\WE\cap\COF$ are {\em acyclic cofibrations};

If the pair of morphisms $i:A\to B, p:X\to Y$ satisfies the condition
(CM 4) we say that $i$ satisfies the {\em left lifting property} (LLP)
with respect to $p$ or that  $p$ satisfies the {\em right lifting 
property} (RLP) with respect to $i$.

\subsection{}
\label{first}

Fix a base ring $k$ and let $C(k)$ be the category of unbounded
complexes over $k$.

Let $\CC$ be a category endowed with a couple of adjoint functors
$$ \#:\CC\rlarrows C(k): F$$
so that $F$ is left adjoint to $\#$.

Suppose that

(H0) $\CC$ admits finite limits and arbitrary colimits;
the functor $\#$ commutes with filtered colimits.

(H1) Let  $d\in\Bbb Z$ and let $M\in C(k)$ be the complex
$$ \cdots\lra 0\lra k = k\lra 0\lra\cdots$$
concentrated in the degrees $d, d+1$. The canonical map
$A\lra A\coprod F(M)$ induces a quasi-isomorphism 
$A^{\#}\lra (A\coprod F(M))^{\#}.$

We define the three classes of morphisms in $\CC$
as follows:

--- $f\in\Mor(\CC)$ belongs to $\WE$ if $f^{\#}$ is a quasi-isomorphism;

--- $f\in\Mor(\CC)$ belongs to $\FIB$ if $f^{\#}$ is (componentwise)
surjective;

--- $f\in\Mor(\CC)$ belongs to $\COF$ if it satisfies the LLP 
with respect to all acyclic fibrations.

\subsubsection{}
\begin{thm}{cmc-structure}
Let a category $\CC$ be endowed with a couple of adjoint functors
$$ F :\CC\rlarrows C(k): \#$$
so that the conditions (H0),(H1) are fulfilled. Then the classes
$\WE,\FIB,\COF$ of morphisms in $\CC$ described above define on $\CC$
a  CMC structure.
\end{thm}

The proof of \Thm{cmc-structure} will be given in~\ref{pf-cmc-str}.

\subsubsection{Adding a variable to kill a cycle}
\label{attaching-cmc}
Let $A\in\CC$,
$M\in C(k)$ and let $\alpha:M\ra A^{\#}$ be a map in $C(k)$ (in particular,
$\alpha$ commutes with the differentials). 

Define a functor 
$$ h_{A,\alpha}:\CC\lra\Ens$$
by the formula
$$ h_{A,\alpha}(B)=\{(f,t)|f:A\ra B\in\Mor(\CC), t\in\Hom^{-1}(M,B^{\#}):
d(t)=f^{\#}\circ\alpha\}.$$

Since $\CC$ admits colimits, the functor $h_{A,\alpha}$ is 
represented as follows. Put $C=\Cone(\alpha)$. One has a couple of maps
$f:A^{\#}\ra C$ and $t\in\CHom^{-1}(M,C)$ satisfying 
$d(t)=f^{\#}\circ\alpha$.
Let now $B$ be a colimit of the diagram
$$ A\lla F(A^{\#})\lra F(C).$$
One sees immediately that the couple of maps $A\ra B$, 
$M\ra F(C)^{\#}\ra B^{\#}$ represents the functor $h_{A,\alpha}$.

The object of $\CC$ representing $h_{A,\alpha}$, will be denoted by
$A\langle M,\alpha\rangle$.

When $M=k[n]$ and $\alpha:M\ra A^{\#}$ takes the generator of $M$ to a cycle
$a\in Z^n(A)$, the representing object is obtained by ``adding a variable
to kill the cycle $a\in Z^n(A^{\#})$''. In this case we will write
$A\langle T; dT=a\rangle$ instead of $A\langle M,\alpha\rangle$.

\subsubsection{Standard cofibrations and standard acyclic cofibrations}

\label{standard:c,ac} 
Let $M$ be a complex of free $k$-modules with zero differential.
For any $A\in\CC$ and any map $\alpha:M\to A^{\#}$ the map

\begin{equation}
A\to A\langle M,\alpha\rangle
\label{sc}
\end{equation}
is a cofibration. 

\begin{defn}{}
A map $A\to B$ is called {\em a standard cofibration} if it is a direct
limit of a sequence
$$ A=A_0\to A_1\to\ldots\to B$$
where each map $A_i\to A_{i+1}$ is as in~(\ref{sc}).
\end{defn}

Let $M$ be a contractible complex of free $k$-modules. Then

\begin{equation}
A\to A\langle M,\alpha\rangle
\label{sac}
\end{equation}
is an acyclic cofibration. 

\begin{defn}{}
A map $A\to B$ is called {\em a standard acyclic cofibration} if it 
is a direct limit of a sequence
$$ A=A_0\to A_1\to\ldots\to B$$
where each map $A_i\to A_{i+1}$ is as in~(\ref{sac}).
\end{defn}

\subsubsection{The proof of~\Thm{cmc-structure}}
\label{pf-cmc-str}
The axioms (CM 1)--(CM 3) are obvious. Also (CM 4)(ii) is immediate.
Let us check (CM 5)(i).

Let $f:A\ra B\in\Mor(\CC)$. For each $b\in B^{\#}$ define
$C_b=A\langle T_b,S_b; dT_b=S_b\rangle$ and let the map
$g_b:C_b\ra B$ be defined by the conditions
$$ g_b^{\#}(T_b)=b;\ g_b^{\#}(S_b)=db.$$ 
Put $C$ to be the coproduct of $C_b$ under $A$ and let $g:C\ra B$
be the corresponding morphism. The map $A\ra C$ is a standard acyclic
cofibration and $g^{\#}$ is surjective.

Now, let us check (CM 5)(ii). 
For this we will construct for a given map $f:A\ra B$ a sequence
$$ A\ra C_0\ra\ldots\ra C_i\ra C_{i+1}\ra\ldots\ra B$$
of standard cofibrations such that

(1) the maps $g_i^{\#}: C_i^{\#}\ra B^{\#}$ are surjective;

(2) $Z(g_i^{\#}): ZC_i^{\#}\ra ZB^{\#}$ are surjective as well
($Z$ denotes the set of cycles);

(3) if $z\in ZC_i^{\#}$ and $g_i^{\#}(z)$ is a boundary in $B^{\#}$
then the image of $z$ in $C_{i+1}^{\#}$ is a boundary.

Then if we put $C=\dirlim C_i$ and $g: C\ra B$ is defined by $g_i$,
the map $A\ra C$ is a cofibration and $g^{\#}$ is clearly a surjective
quasi-isomorphism since the forgetful functor commutes with filtered
colimits.

The object $C_0$ is constructed exactly as in the proof of (CM 5)(i):
one has to join a pair $(T_b,S_b)$ for each element $b\in B^{\#}$ and 
after that to join a cycle corresponding to each cycle in $B^{\#}$.

In order to get $C_{i+1}$ from $C_i$ one has to join to $C_i$ a
variable $T$ for each pair $(z,u)$ with $z\in ZC_i^{\#}, u\in B^{\#}$,
such that $g_i^{\#}(z)=du$. One has to put $dT=z,\ g_{i+1}^{\#}(T)=u$.

Let us prove now (CM 4)(i). The proof of the property (CM 5)(i) implies 
that if $f:A\ra B$ is a weak equivalence then there exists a decomposition
$f=pi$ where $p$ is an acyclic fibration and $i$ is a standard
acyclic cofibration. If $f$ is also a cofibration then according to 
(CM 4)(ii) there exists $j:B\ra C$ making the diagram
\begin{center}
\begin{picture}(3,3)
   \put(0,0){\makebox(1,1){$B$}}
   \put(2,0){\makebox(1,1){$B$}}
   \put(0,2){\makebox(1,1){$A$}}
   \put(2,2){\makebox(1,1){$C$}}
   \put(0.5,2){\vector(0,-1){1}}
   \put(2.5,2){\vector(0,-1){1}}
   \put(1,2.5){\vector(1,0){1}}
   \put(1,0.45){\line(1,0){1}}
   \put(1,0.55){\line(1,0){1}}

   \multiput(0.8,0.8)(0.1,0.1){14}{\makebox(0.1,0.1){$\cdot$}}
   \put(2.2,2.2){\vector(1,1){0}}

   \put(1,2.5){\makebox(1,0.5){$\scriptstyle i$}}
   \put(0,1){\makebox(0.5,1){$\scriptstyle f$}}
   \put(2.5,1){\makebox(0.5,1){$\scriptstyle p$}}
   \put(1,1){\makebox(1.15,0.65){$\scriptstyle j$}}

\end{picture}
\end{center}

commutative.
This proves that any acyclic cofibration is a retract of a standard one
and this immediately implies (CM 4)(i).

Theorem is proven.

Note that the proof is essentially the one given 
in~\cite{q2} for DG Lie algebras or in~\cite{bog} for commutative
DG algebras.

\subsubsection{}
\begin{rem}{retracts} The proof of the Theorem implies the following:

Any acyclic cofibration is a retract of a standard acyclic cofibration. 

Any cofibration is a retract of a standard cofibration.
\end{rem}

\section{Differential homological algebra}
\label{dha-1}
In this Section $k$ is a fixed commutative base ring.

The first application of~\Thm{cmc-structure} provides a CMC structure
on the category of modules $\Mod(A)$ over a dg $k$-algebra $A$.

The constructions of this Section will be generalized in 
\Sect{section-opalg}
to the category of algebras over any $k$-operad. However, even the case
$A=k$ is not absolutely well-known: it provides the category $C(k)$
of unbounded complexes over $k$ with a CMC structure. 

The category $\mod(A)$ admits another, somewhat dual CMC structure.
These two structures are closely related to a homological algebra developed
in~\cite{sp} for the category of sheaves (of modules over a sheaf of
commutative rings) and in~\cite{a} for the category of modules over a dg 
algebra. Another description of the results of this Section can be found
in~\cite{belu}.

\subsection{Models}
\label{dha1-models}

The obvious forgetful functor $\#:\Mod(A)\ra C(k)$ admits the left adjoint
$F=A\otimes_k$. All limits exist in $\Mod(A)$. 
Let $\alpha:M\ra X^{\#}$ be a morphism in $C(k)$. Then
the morphism $\alpha': A\otimes_kM\ra X$ is defined and we have
$$X\langle M,\alpha\rangle=\cone(\alpha').$$
The condition (H1) is trivially fulfilled.

Cofibrant objects in $\Mod(A)$  are exactly direct summands
of semi-free $A$-modules, see~\cite{a}.

The homotopy category of $\Mod(A)$ will be denoted by $D(A)$. Parallelly,
the category $\Mod^r(A)$ of right dg $A$-modules admits the same CMC structure
and the corresponding homotopy category will be denoted by $D^r(A)$.

Note that $D(A)$ is also triangulated, the shift functor and the exact 
triangles being defined in a standard way.

In the special case $A=k$ cofibrant objects of $\Mod(A)$ are exactly 
{\em $\CK$-projective complexes} of Spaltenstein --- see~\cite{sp}. 

\subsection{Tensor product}
The functor
$$ \otimes_A:\Mod^r(A)\times\Mod(A)\ra\Mod(k)$$
is defined as usual: for $M,N\in\Mod(A)$ the tensor product $M\otimes_AN$
is the colimit of the diagram
$$ M\otimes A\otimes N\rra M\otimes N$$
where $\otimes=\otimes_k$, and the arrows take $m\otimes a\otimes n$
to $ma\otimes n$ and to $m\otimes an$ respectively.
Since $\otimes_A$ takes homotopy equivalences to homotopy equivalences, 
and a quasi-isomorphism of cofibrant objects is a homotopy equivalence,
it admits a left derived functor
$$ \otimes_A^{\Left}:D^r(A)\times D(A)\ra D(k).$$
This is the functor defined actually in~\cite{a}. It can be calculated 
using semi-free resolutions with respect to either of the arguments.

\subsection{Base change}
\label{base-dha1}
Let now $f:A\ra B$ be a morphism of dg $k$-algebras. There is a pair of
adjoint functors
$$ f^*:\Mod(A)\rlarrows\Mod(B):f_*$$
where $f_*$ is just the forgetful functor and $f^*(M)=B\otimes_AM$. Since
the functor $f_*$ is exact and $f^*$ preserves cofibrations, one has a
pair of adjoint functors
$$ {\Left}f^*:D(A)\rlarrows D(B):f_*={\Right}f_*.$$

Note that the functor ${\Left}f^*$ commutes with $\otimes^{\Left}$.

\subsubsection{}
\begin{thm}{asso-comparison}
Let $f:A\ra B$ be a quasi-isomorphism of dg algebras. Then the functors
$({\Left}f^*,f_*)$ establish an equivalence of the derived categories $D(A)$
and $D(B)$.
\end{thm}
\begin{pf}
According to~\cite{q1}, \S4, thm.~3, we have to check that if $M$ is a 
cofibrant $A$-module then the map $M\ra f_*f^*(M)$ is a quasi-isomorphism. 
Since any
cofibrant module is a direct summand of a semi-free module, and the functor
$f^*$ commutes with taking cones (i.e., for any $\alpha:M\ra N$
$\Cone(f^*\alpha)=f^*(\Cone(\alpha))$) the result 
immediately follows.

\end{pf}

As an immediate consequence of~\Thm{asso-comparison} we get the following
comparison result which we firstly knew from L.~Avramov (see~\cite{hla},
thm. 3.6.7)
\subsubsection{}
\begin{cor}{}
Let $f:A\ra A'$ be a quasi-isomorphism of dg algebras, $g:M\ra f_*(M')$ and
$h:N\ra f_*(N')$ be quasi-isomorphisms in $\Mod(A)$. Then the induced map
$M\otimes^{\Left}_AN\ra M'\otimes^{\Left}_{A'}N'$ is a quasi-isomorphism.
\end{cor}

\Thm{asso-comparison} will be generalized in \Sect{section-opalg}
to the case of operad algebras --- see \Thm{equivalence}.

\subsection{Flat modules}
Let $A$ be a DG algebra over $k$, $M$ be a $A$-module. We will call $M$
flat (in the terminology of~\cite{a} - $\pi$-flat) if the functor 
$\otimes_AM$ carries quasi-isomorphisms into quasi-isomorphisms. 

\subsubsection{}
\begin{lem}{cofibrant->flat}
1. Any cofibrant $A$-module is $A$-flat.

2. A filtered colimit of flat $A$-modules is $A$-flat.

3. Let $f:X\to Y$  be  a map of flat $A$-modules.  Then the cone
$\Cone(f)$ is also flat.
\end{lem}
\begin{pf}
For the claims 1,2 see~\cite{a}, 6..1, 6.2, 6.6.  The tensor product 
commutes with taking cone --- this implies the third claim.
\end{pf}

\subsubsection{}
\begin{lem}{flatOK}
Let $\alpha:M\ra M'$ be a quasi-isomorphism of flat $A$-modules. Then
for each $N\in\Mod^r(A)$ the map 
$1\otimes\alpha: N\otimes_AM\ra N\otimes_AM'$
is a quasi-isomorphism.
\end{lem}
\begin{pf}
See~\cite{a}, 6.8.
\end{pf}

Thus, the functor $\otimes_A^{\Left}$ can be calculated using flat 
resolutions.

\section{Algebras over an operad}
\label{section-opalg}

\subsection{Introduction}
In this Section we define, using~\Thm{cmc-structure}, a CMC structure
on the category $\Alg(\CO)$ of algebras over an operad $\CO$ which is
{\em $\Sigma$-split} (see~\Defn{ssplit} below). The base 
tensor category is always the category of complexes $C(k)$ over a fixed
commutative ring $k$. All necessary definitions can be found in~\cite{hla},
\S 2,3.

Recall that the forgetful functor $\#:\Alg(\CO)\ra C(k)$ admits a left
adjoint {\em free $\CO$-algebra functor}
$F: C(k)\ra\Alg(\CO)$ which takes a complex $V$ to the $\CO$-algebra
$$FV=\bigoplus_{n\geq 0}(\CO(n)\otimes V^{\otimes n})/\Sigma_n,$$
$\Sigma_n$ being the symmetric group.

\subsubsection{}
\begin{thm}{opalg}
Let $\CO$ be a $\Sigma$-split operad in $C(k)$. The category $\Alg(\CO)$ 
endowed with the couple of adjoint functors
$$ \#:\Alg(\CO)\rlarrows C(k): F$$
satisfies the conditions (H0), (H1). Thus, $\Alg(\CO)$ admits a CMC
structure in which $f:A\ra B$ is a weak equivalence if $f^{\#}$ is a 
quasi-isomorphism and is a fibration if $f^{\#}$ is surjective.
\end{thm}

The proof of~\Thm{opalg} will be given in~\ref{pf-opalg}.

\subsection{$\Sigma$-split operads}

In this subsection we define a class of operads for which \Thm{opalg} 
is applicable.
Let us just mention two important examples of a $\Sigma$-split operad:

--- Any operad in $C(k)$ is $\Sigma$-split if $k\supseteq\Bbb Q$.

--- The operad $\Ass_k$ of associative $k$-algebras is $\Sigma$-split 
    for any $k$.

\subsubsection{Asymmetric operads} We will call an asymmetric operad in $C(k)$
 ``an operad without the symmetric group'': 
it consists of a collection $\CT(n)\in C(k),
\ n\geq 0$, of unit $1:k\ra\CT(1)$, of associative multiplications, but with
no symmetric group action required.

There is a couple of adjoint functors
$$ \Sigma:\asop(k)\rlarrows\op(k):\#$$
between the category of asymmetric operads in $C(k)$ and that of operads.
Here $\#$ is the forgetful functor and $\Sigma$ is defined as follows.

Let $\CT$ be an asymmetric operad. We define 
$\CT^{\Sigma}(n)=\CT(n)\otimes k\Sigma_n$;
multiplication is defined uniquely by the multiplication in $\CT$ in order to
be $\Sigma$-invariant.

For an operad $\CO$ the adjunction map $\pi:\CO^{\#\Sigma}\ra\CO$ is given
by the obvious formula
\begin{equation}
 \pi(u\otimes\sigma)=u\sigma
\label{adj-pi}
\end{equation}
where $u\in\CO(n), \sigma\in\Sigma_n$.

\subsubsection{Notations: symmetric groups}
In this subsection we denote by $\langle n\rangle$ the ordered set 
$\{1,\ldots,n\}$. Let $f:\langle s\rangle\lra\langle n\rangle$ be an
injective monotone map. This defines a monomorphism $\iota_f:
\Sigma_s\ra\Sigma_n$ in the obvious way: for $\rho\in\Sigma_s$
\begin{equation}
\iota_f(\rho)(i)=\begin{cases}
                  i\text{ if } i\not\in f(\langle s\rangle)\\
                  f(\rho(j))\text{ if } i=f(j) 
                 \end{cases}
\end{equation}
Define a map (not  a homomorphism) $\rho_f:\Sigma_n\ra\Sigma_s$ by the 
condition
\begin{equation}
\rho=\rho_f(\sigma)\text{ iff } \rho(i)<\rho(j)\Leftrightarrow \sigma(f(i))<\sigma(f(j))
\end{equation}
Define a set $T_f\subseteq\Sigma_n$ by 
$$ T_f=\{\sigma\in\Sigma_n|\sigma\circ f:\langle s\rangle\ra\langle n\rangle
\text{ is monotone }\}. $$
\begin{lem}{}
For $\sigma\in\Sigma_n$ there is a unique presentation
$$ \sigma=\tau\iota_f(\rho)$$
with $\tau\in T_f$ and $\rho\in\Sigma_s$.
In this presentation $\rho=\rho_f(\sigma)$.
\end{lem}
\begin{pf}
Obvious.
\end{pf}

For $M\in C(k)$ we will write $M\otimes\Sigma_n$ instead of 
$M\otimes_kk\Sigma_n$. Also if $M$ is a right $\Sigma_n$-module, $N$ is a left 
$\Sigma_n$-module and $\Sigma$ is a subgroup in $\Sigma_n$ then we write
$M\otimes_{\Sigma}N$ instead of $M\otimes_{k\Sigma}N=(M\otimes N)/\Sigma.$

If $M$ admits a right $\Sigma_n$-action and
$f:\langle s\rangle\ra\langle n\rangle$ is monotone injective, the
map $\rho_f$ defined above induces a map 
$$ M\otimes\Sigma_n\lra M\otimes\Sigma_s$$
carrying the element $m\otimes\sigma$ to $m\tau\otimes\rho_f(\sigma)$
where, as in Lemma above, $\sigma=\tau\iota_f(\rho_f(\sigma))$. 
We denote this map also by $\rho_f$.

Note that the map $\rho_f$ is equivariant with respect to right 
$\Sigma_s$-action.

\subsubsection{Notations: operads}
An operad $\CO$ is defined by a collection of multiplication maps
$$\gamma:\CO(n)\otimes\CO(m_1)\otimes\cdots\otimes\CO(m_n)\lra
\CO(\sum m_i).$$
One defines the operations $\circ_k:\CO(n)\otimes\CO(m)\ra\CO(m+n-1)$
to be the compositions
$$ \CO(n)\otimes\CO(m)\ra\CO(n)\otimes\CO(1)^{\otimes k-1}\otimes\CO(m)
\otimes\CO(1)^{\otimes n-k}\overset{\gamma}{\ra}\CO(m+n-1)$$
with the first map induced by the units $1:k\ra\CO(1)$.
The multiplications $\circ_k$ can be described one through another
for different $k$ using the symmetric group action on $\CO(n)$.

\subsubsection{$\Sigma$-split operads}$ $

$\Sigma$-splitting of an operad $\CO$ is a collection of maps
$\CO(n)\to\CO^{\#\Sigma}(n)$ which splits the adjunction map
$\pi:\CO^{\#\Sigma}\to\CO$ from~(\ref{adj-pi}). Of course, there is some
condition describing a compatibility of these maps with the multiplications
maps. Here is the definition.

\begin{defn}{ssplit}
(1) Let $\CO$ be an operad in $C(k)$. {\em
$\Sigma$-splitting } of $\CO$ is a collection of maps of complexes 
$t(n):\CO(n)\ra\CO^{\#\Sigma}(n)$ such that 

(EQU) $t(n)$ is $\Sigma_n$-equivariant

(SPL) $\pi\circ t(n)=\id:\CO(n)\ra\CO(n)$ and

(COM) for any $m,n>0$ and $1\leq k\leq n$ the diagram
\begin{center}
\begin{picture}(13,7)
   \put(0,0){\makebox(3,1){$\CO(n)\otimes\CO(m)\otimes\Sigma_{n-1}$}}
   \put(0,6){\makebox(3,1){$\CO(n)\otimes\CO(m)$}}
   \put(0,4){\makebox(3,1){$\CO(n)\otimes\Sigma_n\otimes\CO(m)$}} 
   \put(0,2){\makebox(3,1){$\CO(n)\otimes\Sigma_{n-1}\otimes\CO(m)$}}  
   \put(7,0){\makebox(3,1){$\CO(n+m-1)\otimes\Sigma_{n-1}$}}
   \put(7,3){\makebox(3,1){$\CO(n+m-1)\otimes\Sigma_{n+m-1}$}} 
   \put(7,6){\makebox(3,1){$\CO(n+m-1)$}}

   \put(4,0.5){\vector(1,0){2}}
   \put(3.5,6.5){\vector(1,0){3}}

   \put(1.5,6){\vector(0,-1){1}}
   \put(1.5,4){\vector(0,-1){1}}  
   \put(1.5,2){\vector(0,-1){1}}  
   \put(8.5,6){\vector(0,-1){2}}
   \put(8.5,3){\vector(0,-1){2}}   

   \put(3,6.5){\makebox(4,1){$\scriptstyle \circ_k$}}
   \put(3,0.5){\makebox(4,1){$\scriptstyle \circ_k\otimes\id$}}

   \put(1.5,5){\makebox(2,1){$\scriptstyle t\otimes\id$}}
   \put(1.5,3){\makebox(2,1){$\scriptstyle \rho_f\otimes\id$}}     
   \put(1.5,1){\makebox(2,1){$\scriptstyle \sigma_{23}$}}        
   \put(8.5,4){\makebox(2,2){$\scriptstyle t$}}      
   \put(8.5,1){\makebox(2,2){$\scriptstyle\rho_g$}}
\end{picture}
\end{center}

is commutative.

Here $f:\langle n-1\rangle\ra\langle n\rangle$ is the map omitting the
value $k$, $g:\langle n-1\rangle\ra\langle m+n-1\rangle$ 
omits the values $k,\ldots,k+m-1$ and $\sigma_{23}$ is the standard twist
interchanging the second and the third factors.

(2) An operad $\CO$ is $\Sigma$-split if it admits a $\Sigma$-splitting.
\end{defn}

\begin{rem}{}  It is sufficient to require the validity of (COM) only for,
say, $k=1$. The compatibility of the map $t$ with other multiplications
$\circ_k$ then follows immediately since $\circ_k$ can be expressed through
$\circ_1$ using the symmetric group action on the components of the operad.
\end{rem}

\subsubsection{Examples}
There are two very important examples of $\Sigma$-split operads.

1. Let $\CT$ be an asymmetric operad and $\CO=\CT^{\Sigma}$.
Then the composition
$$ \CO(n)=\CT^{\Sigma}(n)\lra\CO^{\#\Sigma}(n)$$
defines a $\Sigma$-splitting of $\CO$. 

Let $\Com_k$ be the operad given by 
$\Com_k(n)=k$ for all $n$.The action of $\Sigma_n$ on $k$ is supposed to be
trivial. The algebras over $\Com_k$ are just commutative dg $k$-algebras:
$\Alg(\Com_k)=\cdga(k)$. Put $\Ass_k=\Com_k^{\#\Sigma}$. One has
$\Alg(\Ass_k)=\dga(k)$, the category of associative dg algebras,
and $\Ass_k$ is naturally $\Sigma$-split. 

2. Suppose $k\supseteq\Bbb Q$. Then any operad in $C(k)$ is $\Sigma$-split: 
the splitting is defined by the formula
$$t(u)=\frac{1}{n!}\sum_{\sigma\in\Sigma_n}u\sigma^{-1}\otimes\sigma.$$

The operad $\Com_k$ is $\Sigma$-split only when $k\supseteq\Bbb Q$.
The same is true for the operad $\Lie_k$ such that the algebras over $\Lie_k$
are just dg Lie $k$-algebras. We denote in the sequel by $\dgl(k)=\Alg(\Lie_k)$
the category of dg Lie algebras over $k$. 

\subsection{Extension of a homotopy on free algebras}
\label{hfree}

Let $V\in C(k)$. Let $\alpha:V\ra V$ be a map of complexes of degree zero
and let $h\in\Hom(V,V[-1])$ satisfy the property
$$ d(h)=\id_V-\alpha.$$
The endomorphism $\alpha$ induces the endomorphism $F(\alpha)$ of the free
$\CO$-algebra $F(V)$ by the obvious formula
$$ F(\alpha)(u\otimes x_1\otimes\cdots\otimes x_n)=u\otimes\alpha(x_1)
\otimes\cdots\otimes\alpha(x_n).$$

We will describe now a nice homotopy $H$ connecting $\id_{F(V)}$ with 
$F(\alpha)$. This is a sort of "skew derivation" on $F(V)$ defined by $h$.

The restriction of $H$ on $F_n(V)=\CO(n)\otimes_{\Sigma_n}V^{\otimes n}$
is given by the composition
\begin{equation}
\CO(n)\otimes_{\Sigma_n}V^{\otimes n}
\overset{t}{\ra}\CO(n)\otimes\Sigma_n\otimes_{\Sigma_n}V^{\otimes n}=
\CO(n)\otimes V^{\otimes n}@>\sum_{p+q=n-1}\id\otimes\alpha
^{\otimes p}\otimes h\otimes\id^{\otimes q}>>
\CO(n)\otimes_{\Sigma_n}V^{\otimes n}.
\end{equation}

The property $dH=\id_{F(V)}-F(\alpha)$ is verified immediately.

\subsection{Proof of~\Thm{opalg}}
\label{pf-opalg}
The property (H0) is obvious. 

Let us prove (H1). 
Let $A$ be a $\CO$-algebra and let $X\in C(k)$ be a contractible complex.
Put $V=A^{\#}\oplus X$. The complex $V$ is homotopy equivalent to $A^{\#}$,
the maps between $A^{\#}$ and $V$ being obvious and the homotopy equivalence 
being defined by a map $h:V\ra V$ of degree $-1$ which vanishes on $A^{\#}$.  
One has
$dh=\alpha$ where $\alpha:V\ra V$ is the composition 
$$V=A^{\#}\oplus X\ra A^{\#}\ra V.$$

According to~\ref{hfree}, $h$ defines a homotopy
$$H: F(V)\lra F(V)$$
of degree $-1$ extending $h$.

Let now  $I$ be the kernel of the natural projection
$F(A^{\#})\ra A$. Let $J$ be the ideal in $F(V)$ generated by $I$. 
We will prove now that the homotopy $H$ satisfies the property
\begin{equation}
H(J)\subseteq J.
\label{J-H-invariant}
\end{equation}
Then $H$ induces a homotopy on $F(V)/J=A\coprod F(X)$ which proves the
theorem.

To prove the property~(\ref{J-H-invariant}) let us consider the restriction of
$H$ to $\CO(n)\otimes A^{\otimes r}\otimes V^{\otimes s}$ with $n=r+s$.

An easy calculation using the properties $\alpha|_A=\id,\ h|_A=0$ shows
that this restriction of $H$ can be calculated as the composition
\begin{multline}
\CO(n)\otimes A^{\otimes r}\otimes V^{\otimes s}\ra
\CO(n)\otimes\Sigma_n\otimes A^{\otimes r}\otimes V^{\otimes s}
\overset{\rho}{\ra}\\
\CO(n)\otimes\Sigma_s\otimes A^{\otimes r}\otimes V^{\otimes s}
\overset{\tau_{23}}{\ra}\CO(n)\otimes A^{\otimes r}\otimes\Sigma_s
\otimes V^{\otimes s}\ra \\
\CO(n)\otimes A^{\otimes r}\otimes V^{\otimes s}@>\sum_{p+q=s-1}
\id\otimes\id\otimes\alpha^{\otimes p}\otimes h\otimes\id^{\otimes q}>>
\CO(n)\otimes A^{\otimes r}\otimes V^{\otimes s}\ra
\CO(n)\otimes_{\Sigma_n}V^{\otimes n}.
\label{Hav}
\end{multline}

To check~(\ref{J-H-invariant}) note that the ideal $I\subseteq F(A^{\#})$ 
is generated over $k$ by the expressions
$$b\otimes x_1\otimes\cdots\otimes x_m-
\mu(b\otimes x_1\otimes\cdots\otimes x_m)$$
with $b\in\CO(m), x_i\in A$. Therefore  the ideal $J\subseteq F(V)$
is generated over $k$ by the expressions
\begin{equation}
a\circ_1 b\otimes x_1\cdots\otimes x_m\otimes y_1\otimes\cdots\otimes y_{n-1}
-a\otimes\mu(b\otimes x_1\otimes\cdots\otimes x_m)\otimes 
y_1\otimes\cdots\otimes y_{n-1}.
\label{genJ}
\end{equation}
Hence, we have to check that $H$ transforms an element of form~(\ref{genJ})
into an element of $J$. This easily follows from the axiom (COM) and
formula~(\ref{Hav}).

Theorem is proven.

\subsection{Notations} The homotopy category of $\Alg(\CO)$ is denoted by
$\hoalg(\CO)$. For the special values of $\CO$ we denote
$\hoass(k)=\hoa(\Ass_k),\ \hocom(k)=\hoa(\Com_k),\ \hol(k)=\hoa(\Lie_k)$.

\subsection{Base change}
Consider now a map $\alpha: \CO\ra \CO'$ of operads. We will study
direct and inverse image functors between the categories $\Alg(\CO)$ and
$\Alg(\CO')$.

This generalizes the considerations of~\ref{asso-comparison} to the case of 
operad algebras.

\subsubsection{Direct image}

Let $A$ be a $\CO'$-algebra. Its direct image $\alpha_*(A)$ is just the 
$\CO$-algebra obtained from $A$ by forgetting ``the part of structure'':
the multiplication map is given by the composition
$$ \CO(n)\otimes A^{\otimes n}\overset{\alpha\otimes 1}{\lra}
 \CO'(n)\otimes A^{\otimes n}\lra A.$$

This functor is obviously exact.

\subsubsection{Inverse image}

The inverse image $\alpha^*:\Alg(\CO)\ra\Alg(\CO')$ is by definition the
functor left adjoint to $\alpha_*$.

Let us explicitly construct $\alpha^*$. Let $F$ and $F'$ be the free
$\CO$-algebra and free $\CO'$-algebra functors respectively.
For $A\in\Mod(\CO)$ let $I_A$ be the kernel of the natural map
$F(A^{\#})\ra A$. Then, if $F(\alpha): F(A^{\#})\ra F'(A^{\#})$ is the map 
induced by $\alpha$, one defines
$$ \alpha^*(A)=F'(A^{\#})/(F(\alpha)(I_A)).$$ 

\subsubsection{Derived functors} 

We wish now to construct an adjoint pair of derived functors
$$ \Left\alpha^*:\hoalg(\CO)\rlarrows\hoalg(\CO'):
\Right\alpha_*=\alpha_*.$$
Let us check the conditions of~\cite{q1}, \S 4, thm.~3.

Let $M\in C(k), A\in\Alg(\CO)$. Any map $f:M\ra A^{\#}$  defines a map
$f': M\ra(\alpha^*A)^{\#}$. Then one immediately sees that there is
a canonical isomorphism
$$ \alpha^*(A\langle M,f\rangle)=\alpha^*(A)\langle M,f\rangle$$
since these two $\CO'$-algebras just represent isomorphic functors.
This immediately implies that $\alpha^*$ carries 
standard cofibrations to standard cofibration and standard acyclic cofibrations
to standard acyclic cofibrations. 
Then~\Rem{retracts} implies that $\alpha^*$  preserves cofibrations and
acyclic cofibrations. 

Let us check that $\alpha^*$ carries fibrations to fibrations. In fact,
if $f:A\ra B$ is a fibration, then $f^{\#}: A^{\#}\ra B^{\#}$ is surjective.
Thus the induced map $F'(A^{\#})^{\#}\ra F'(B^{\#})^{\#}$ is also surjective
which ensures that the induced map of the quotients 
$$ \alpha^*(A)^{\#}\ra \alpha^*(B)^{\#}$$
is also surjective.

Let us check that $\alpha^*$ carries homotopy equivalences to weak 
equivalences.

In fact, if in the diagram
$$ X\overset{i}{\lra}X^I\overset{p}{\lra}X\times X$$
$i$ is a trivial cofibration and $p$ is a fibration, the functor
$\alpha^*$ gives rise to the  diagram
$$\alpha^*(X)\overset{\alpha^*(i)}{\lra}\alpha^*(X^I)
\overset{\alpha^*(p)}{\lra}\alpha^*(X\times X)\overset{j}{\lra}
\alpha^*(X)\times\alpha^*(X).$$

The map $\alpha^*(i)$ is already known to be acyclic cofibration.
Thus, any pair of maps from somewhere to $\alpha^*(X)$ defined by a map 
to $\alpha^*(X^I)$, induce the same map in homology. This implies
that if $f:X\ra Y$ is a homotopy equivalence in $\Alg(\CO)$ then
$\alpha^*(f)$ is an isomorphism in $\hoalg(\CO')$.

Since weak equivalences of cofibrant objects in $\Alg(\CO)$ are homotopy
equivalences so $\alpha^*$ carries them to weak equivalences in $\Alg(\CO').$

This proves the following
\subsubsection{}
\begin{thm}{derived^*}
Inverse and direct image functors define a pair of adjoint derived functors
$$ \Left\alpha^*:\hoalg(\CO)\rlarrows\hoalg(\CO'):
\Right\alpha_*=\alpha_*.$$
\end{thm}

\subsection{Equivalence} 
\label{eq-opalg}
Suppose that $\alpha:{\CO}\ra{\CO'}$ is a quasi-isomorphism of $\Sigma$-split
operads {\em compatible with the splittings}. We shall prove
that $\Left\alpha^*$ and $\alpha_*$ establish an equivalence of the
homotopy categories $\hoalg(\CO)$ and $\hoalg(\CO')$.

In order to do this, one has to check that the adjunction map
$$ \eta_A: A\lra \alpha_*(\alpha^*(A))$$
is a weak equivalence for any cofibrant $A$.

\subsubsection{1st reduction} Since  retract of a weak equivalence is a weak
equivalence, it suffices to prove the assertion for $A$ standard cofibrant.

\subsubsection{2nd reduction} Since any standard cofibrant object is a filtered
colimit of finitely generated ones, and the functors $\alpha^*$ and $\alpha_*$ 
commute with filtered colimits, it suffices to prove that $\eta_A$ is
a weak equivalence when $A$ is a finitely generated standard cofibrant 
algebra.

\subsubsection{}
\label{fg-cof}
Let now $A$ be a finitely generated standard cofibrant algebra. Let
$\{x_i\}_{i\in I}$ be a finite set of (graded free) generators of $A$.
Choose a full order on the set $I$ of generators in order that for any
$i\in I$ the differential $d(x_i)$ belongs to the algebra generated by 
$x_j, j<i$. 

For any multi-index $m:I\ra\Bbb N$ denote $|m|=\sum m_i$. Denote by
$M$ the set of all multi-indices. The set $M$ is well-ordered
with respect to "inverse lexicographic" order:

$$ m>m'\text{ if there exists }i\in I\text{ so that } m_j=m'_j\text{ for }
j>i\text{ and }m_i>m'_i.$$

Then the algebra $A$ as a graded $k$-module is a direct sum indexed by $M$
of the components 
$$\CO(|m|)\otimes_{\Sigma_m}\bigotimes_{i\in I}x_i^{\otimes m_i}.$$
Here $\Sigma_m=\prod_{i\in I}\Sigma_{m_i}$.

This defines an increasing filtration of $A$ 
$$ \CF_d(A)=\sum_{m<d}\CO(|m|)\otimes_{\Sigma_m}
\bigotimes_{i\in I}x_i^{\otimes m_i}$$
indexed by $M$ which is obviously a filtration by subcomplexes.

The functor $\alpha^*:\Alg(\CO)\ra\Alg(\CO')$ commutes with the functor
forgetting the differentials. Thus, $A'=\alpha^*(A)$ admits the filtration
analogous to $\{\CF_d(A)\}_{d\in M}$. 

In order to prove that the map $\eta_A$ is a weak equivalence, we will
prove by induction that the map
$$\CF_d(A)\lra\CF_d(A')$$
is quasi-isomorphism. For this one has to check that the maps
$$ \gr_d(\eta):\gr_d(A)\lra\gr_d(A')$$
are quasi-isomorphisms where
$$ \gr_d(A)=\CF_{d+1}(A)/\CF_d(A)\approx \CO(|d|)\otimes_{\Sigma_m}k$$        
and similarly $\gr_d(A')\approx \CO'(|d|)\otimes_{\Sigma_m}k$. 

Now we will use that the map $\alpha:\CO\ra\CO'$
is compatible with $\Sigma$-splittings. In fact, in this case the map
$\gr_d(\eta)$ is a retract of the map
$$ \CO^{\Sigma}(|d|)\otimes_{\Sigma_d}k\lra\CO^{'\Sigma}(|d|)
\otimes_{\Sigma_d}k$$
which is obviously a quasi-isomorphism.

Thus we have proven the following

\subsubsection{}
\begin{thm}{equivalence}
Let $\alpha:\CO\ra\CO'$ be a quasi-isomorphism of $\Sigma$-split operads
compatible with splittings. Then the functors
$$ \Left\alpha^*:\hoa(\CO)\rlarrows\hoa(\CO'):\Right\alpha_*=\alpha_*$$
are equivalences of the homotopy categories.
\end{thm}

\subsection{Simplicial structure on $\Alg(\CO)$}
From now on the base ring $k$ is supposed to contain the rationals.
We define on $\Alg(\CO)$ the structure of simplicial category which is
a direct generalization of the definitions~\cite{bog}, Ch.~5.

\subsubsection{Polynomial differential forms}

Recall (cf.~\cite{bog},~\cite{hdtc}, ch.~6) the definition of simplicial
commutative dg algebra $\Omega=\{\Omega(n)\}_{n\geq 0}$.

For any $n\geq 0$ the dg algebra $\Omega_n$ is the algebra of polynomial
differential forms on the standard $n$-simplex $\Delta(n)$.

Thus, one has 
$$\Omega_n=k[t_0,\ldots,t_n,dt_0,\ldots,dt_n]/(\sum t_i-1,\sum dt_i).$$

The algebras $\Omega_n$ form a simplicial commutative dg algebra:
a map $u:[p]\to[q]$ induces the map $\Omega(u):\Omega_q\to\Omega_p$
defined by the formula $\Omega(u)(t_i)=\sum_{u(j)=i}t_j$.

\subsubsection{Functional spaces for $\CO$-algebras}
Let $A,B\in\Alg(\CO)$. We define $\Hom^{\Delta}(A,B)\in\Delta^0\Ens$
to be the simplicial set whose $n$ simplices are
$$\Hom^{\Delta}_n(A,B)=\Hom(A,\Omega_n\otimes B).$$
Note that $\Omega_n$ being a commutative dg algebra over $k$, the tensor
product admits a natural $\CO$-algebra structure.

\subsubsection{}
\begin{lem}{l5.2}(cf.~\cite{bog},~Lemma 5.2)
There is a natural morphism
$$\Phi(W):\Hom(A,\Omega(W)\otimes B)\iso\Hom_{\Delta^0\Ens}
(W,\Hom^{\Delta}(A,B))$$
which is a bijection provided $W$ is finite.
\end{lem}
\begin{pf}
The map $\Phi$ is defined in a standard way.
One has obviously that $\Phi(\Delta(n))$ is a bijection for any $n$.
Now, the contravariant functor $\Omega:\Delta^0\Ens\to\cdga(k)$
carries  colimits to  limits; the functor
$$\_\otimes B:\cdga(k)\to\Alg(\CO)$$
preserves finite limits. This proves that $\Phi(W)$ is bijection for
any finite simplicial set $W$.
\end{pf}

\subsubsection{}
\begin{lem}{p5.3}
Let $i:A\to B$ be a cofibration and $p:X\to Y$ be a fibration 
in $\Alg(\CO)$. Then the canonical map
$$ (i^*,p_*):\Hom^{\Delta}(B,X)\to\Hom^{\Delta}(A,X)
\times_{\Hom^{\Delta}(A,Y)}\Hom^{\Delta}(B,Y)$$
is a Kan fibration. It is acyclic if $i$ or $p$ is acyclic.
\end{lem}
\begin{pf}
See the proof of~\cite{bog}, Prop.~5.3.
\end{pf}

The assertions below immediately follow from~\Lem{p5.3}, see also~\cite{bog},
Ch.~5.
\subsubsection{}
\begin{cor}{} 
Let $i:A\to B$ be a cofibration and $C\in\Alg(\CO)$. Then
$$ i^*:\Hom^{\Delta}(B,X)\to\Hom^{\Delta}(A,X)$$
is a Kan fibration. It is acyclic if $i$ is acyclic.
\end{cor}

\subsubsection{}
\begin{cor}{}
If $A$ is cofibrant then $\Hom^{\Delta}(A,X)$ is Kan for every $X$.
\end{cor}

\subsubsection{}
\begin{cor}{}
If $A$ is cofibrant and $p:X\to Y$ is fibrant then
$p_*:\Hom^{\Delta}(A,X)\to\Hom^{\Delta}(A,Y)$ is Kan fibration. It is acyclic
if $p$ is acyclic.
\end{cor}

\subsubsection{}
\begin{cor}{}
Let $A$ be cofibrant and let $f:X\to Y$ be weak equivalence. Then
$$f_*:\Hom^{\Delta}(A,X)\to\Hom^{\Delta}(A,Y)$$
is a weak equivalence.
\end{cor}

\subsubsection{}
\begin{rem}{}Note that the canonical map $A\to A^I$ is not usually a
cofibration:
take, for instance, $A$ to be the trivial (one-dimensional) Lie algebra.
Then $A^I$ is commutative and of course is not cofibrant.
\end{rem}

\subsubsection{Simplicial homotopy}   $ $

\begin{defn}{}Two maps $f,g:A\to B$ in $\Alg(\CO)$ are called {\em simplicially 
homotopic} if there exists $F\in\Hom_1^{\Delta}(A,B)$ such that $d_0F=f,\ d_1F=g$.
\end{defn}

All the assertions of~\cite{bog}, Ch.~6, are valid in our case. In particular,
simplicial homotopy is an equivalence relation provided $A$ is cofibrant. In 
this case simplicial homotopy coincides with both right and left homotopy relations
defined in~\cite{q1}. This allows one to realize the homotopy category $\hoa(\CO)$
as the category having the cofibrant $\CO$-algebras as the objects and
the set $\pi_0\Hom^{\Delta}(A,B)$ as the set of morphisms from $A$ to $B$.

It seems however that the simplicial category $\Alg^{\Delta}(\CO)$ defined by

--- $\Ob\Alg^{\Delta}(\CO)$ is the collection of cofibrant 
$\CO$-algebras;

---  $A,B\mapsto\Hom^{\Delta}(A,B)$

is more useful then the homotopy category $\hoa(\CO)$.

\section{Modules over operad algebras}
\label{moa}
In this Section we study the category of modules over an operad algebra
$(\CO,A)$.
This can be described as the category of modules over the universal 
enveloping algebra  $U(\CO,A)$. The corresponding derived category 
$DU(\CO,A)$ can be different for quasi-isomorphic operad algebras, so
one has to "derive" this construction to get an invariant depending
only on the isomorphism class of $(\CO,A)$ in the homotopy category 
$\hoa(\CO)$. To get this, one should substitute the algebra $A$ with 
its cofibrant resolution --- thus substituting $A$-modules with 
"virtual $A$-modules" and the enveloping algebra of $A$ --- 
with the "derived enveloping algebra".

Starting from~\ref{derived-e-a} we suppose that the operad $\CO$ is
$\Sigma$-split.

\subsection{Modules. Enveloping algebra.} 

We refer to~\cite{hla}, ch.~3, for the definition of $(\CO,A)$-modules,
$(\CO,A)$-tensor algebra $T(\CO,A)$ and the universal enveloping algebra 
$U(\CO,A)$.

\subsubsection{}
\begin{lem}{U-commutes-with-fcolimits}The functor $U(\CO,\_):\Alg(\CO)\ra
\dga(k)$ commutes with filtered colimits.
\end{lem}
\begin{pf}
Recall that the enveloping algebra $U(\CO,A)$ coincides with the colimit
(both in  $\dga(k)$ and in $C(k)$) of the diagram
$$ T(\CO,\#F(\#A))\rra T(\CO,\#A)$$
where $T(\CO,\_)$ is the $\CO$-tensor algebra functor and $F(\_)$ is the free
$\CO$-algebra functor --- see~\cite{hla}, ch.4.
Now the lemma immediately follows from the fact that the functors
$F,T,\#$ commute with filtered colimits.
\end{pf}

\subsection{Functoriality} Let $f=(\alpha,\phi):(\CO,A)\ra(\CO',A')$ be 
a map of operad algebras, where $\alpha:\CO\ra\CO'$ is a map of operads and
$\phi:A\ra\alpha_*(A')$ is a map of $\CO$-algebras. This induces 
a map $U(f):U(\CO,A)\ra U(\CO',A')$ of the corresponding enveloping
algebras and so  by~\ref{base-dha1} one has the following pairs of
adjoint functors
\begin{equation}
f^*:\mod(\CO,A)\rlarrows\mod(\CO',A'): f_*
\label{adj-mod}
\end{equation}

\begin{equation}
{\Left}f^*:DU(\CO,A)\rlarrows DU(\CO',A'): f_*={\Right}f_*.
\label{adj-D}
\end{equation}

The adjoint functors~(\ref{adj-D}) are equivalences provided 
$f:U(\CO,A)\to U(\CO',A')$ is a quasi-isomorphism. Unfortunately, this is
not always the case even when $\alpha$ and $\phi$ are quasi-isomorphisms.

\subsection{Derived enveloping algebra}
\label{derived-e-a} 
Fix a $\Sigma$-split operad $\CO$ in $C(k)$; we will write $U(A)$
instead of $U(\CO,A)$ and $T(V)$ instead of $T(\CO,V)$.

\subsubsection{}
\begin{prop}{UA=UA-X}
Let $A\in\Alg(\CO)$ be cofibrant, $X\in C(k)$ be contractible, 
$A'=A\coprod F(X)$. Then the natural map $U(A)\lra U(A')$ is a 
quasi-isomorphism.
\end{prop}
\begin{pf}
The proof is similar to that of~\ref{eq-opalg}.

{\em 1st reduction. } It suffices to prove the claim for standard cofibrant $A$
since a retract of quasi-isomorphism is quasi-isomorphism.

{\em 2nd reduction.} We can suppose that $A$ is finitely generated since
filtered colimit of quasi-isomorphisms is quasi-isomorphism and the functor $U$
commutes with filtered colimits --- see~\ref{U-commutes-with-fcolimits}.

{\em 3rd step.} (compare with~\ref{eq-opalg}). Let $\{x_i\}_{i\in I}$ be a set
of homogeneous generators of $A$ with $I$ ordered as in~\ref{eq-opalg}.
Let $M$ be the set of multi-indices $m:I\ra{\Bbb N}$ with the 
``opposite-to-lexicographic'' order. Then $U(A)$ as a graded $k$-module
takes form
$$ U(A)=\bigoplus_{m\in M}\CO(|m|+1)\otimes_{\Sigma_m}\bigotimes_ix_i^{\otimes
m_i}$$
where, as in~\ref{eq-opalg}, $\Sigma_m=\prod_{i\in I}\Sigma_{m_i}$.
This defines a filtration of $U(A)$ by subcomplexes
\begin{equation}
\CF_d(U(A))=\sum_{m<d}\CO(|m|+1)\otimes_{\Sigma_m}\bigotimes_ix_i^{\otimes
m_i}.
\label{FU}
\end{equation}

In a similar way, $U(A')$ as a graded $k$-module is isomorphic to a 
tensor algebra; it admits a direct sum decomposition as follows
$$ U(A')=\bigoplus_{m\in M}(\bigoplus_{k\geq 0}\CO(|m|+k+1)\otimes
_{\Sigma_k\times\Sigma_m}X^{\otimes k}\otimes
\bigotimes_ix_i^{\otimes m_i}) .$$

This defines a filtration of $U(A')$ by subcomplexes
$$\CF(U(A'))=\sum_{m<d}(\bigoplus_{k\geq 0}\CO(|m|+k+1)\otimes
_{\Sigma_k\times\Sigma_m}X^{\otimes k}\otimes
\bigotimes_ix_i^{\otimes m_i}) .$$

The associated graded complexes take form
$$ \CO(|d|+1)\otimes_{\Sigma_d}k\text{ for }\gr_d\CF(U(A))$$
and
$$\sum_{k\geq 0}\CO(|d|+k+1)\otimes_{\Sigma_k\times\Sigma_d}X^{\otimes k}
\otimes k\text{ for }\gr_d\CF(U(A')).$$

We have to check that the summands corresponding to $k>0$ are 
contractible. This immediately follows from the contractibility of $X$
and $\Sigma$-splitness of $\CO$.
\end{pf}

\subsubsection{}
\begin{cor}{}
Let $f:A\ra B$ be an acyclic cofibration in $\Alg(\CO)$ with $A$ (and hence 
$B$) cofibrant. Then $U(f)$ is quasi-isomorphism.
\end{cor}
\begin{pf}
Any acyclic cofibration is a retract of a standard one; since everything
commutes with filtered colimits, we immediately get the assertion.
\end{pf}

\subsubsection{}
\begin{cor}{U(w-c)}
Let $f:A\ra B$ be a weak equivalence of cofibrant algebras in $\Alg(\CO)$.
then $U(f)$ is quasi-isomorphism.
\end{cor}
\begin{pf}
Let $B\overset{\alpha}{\lra}B^I\rra B$  be a path diagram for $B$ so that
$\alpha$ is an acyclic cofibration. By Lemma above $U(\alpha)$ is 
quasi-isomorphism. This immediately implies that if $f,g:A\ra B$ are homotopic
then $U(f),U(g)$ induce the same map in cohomology.

Now, if $f:A\ra B$ is a weak equivalence and $A, B$ are cofibrant then $f$
is homotopy equivalence, i.e. there exist $g:B\ra A$ such that the compositions
are homotopic to appropriate identity maps. This implies that $U(f)$ and
$U(g)$ induce mutually inverse maps in the cohomology. In particular, $U(f)$
is quasi-isomorphism.
\end{pf}

\Cor{U(w-c)} allows one to define the left derived functor
$$ {\Left}U:\hoalg(\CO)\to\hoass(k) $$
from the homotopy category of $\CO$-algebras to the homotopy category of
associative dg $k$-algebras.

\subsubsection{}
\begin{lem}{uniqueness-of-d-i}1. Let $f:A\ra B$ be a weak equivalence of 
cofibrant algebras in $\Alg(\CO)$. Then the functors
\begin{equation}
{\Left}f^*:DU(\CO,A)\rlarrows DU(\CO,B): f_*={\Right}f_*.
\end{equation}
of~(\ref{adj-D}) are equivalences.

2. Let $f,g:A\to B$ be homotopic maps of cofibrant algebras. Then there
is an isomorphism of functors 
$$ f_*\iso g_*: DU(\CO,A)\to DU(\CO,B).$$
This isomorphism depends only on the homotopy class of the homotopy
connecting $f$ with $g$.
\end{lem}
\begin{pf}
The first part follows from \Cor{U(w-c)} and~\ref{asso-comparison}.

Let $B\overset{\alpha}{\lra}B^I\overset{p_0,p_1}{\rra} B$  be a path diagram
for $B$ so that $\alpha$ is an acyclic cofibration. Since the functors $p_{0*}$
and $p_{1*}$ are both quasi-inverse to $\alpha_*$, they are naturally 
isomorphic. 
Therefore, any homotopy $F:A\to B^I$ between $f$ and $g$ defines an
isomorphism $\theta_F$ between $f_*$ and $g_*$. Let now
$F_0,F_1:A\to B^I$ be homotopic. The homotopy can be realized by a map
$h:A\to C$ where $C$ is taken from the path diagram
\begin{equation}
B^I\overset{\beta}{\lra}C\overset{q_0\times q_1}{\lra} B^I\times_{B\times B}B^I
\label{2nd-path}
\end{equation}
where $\beta$ is an acyclic cofibration, $q_0\times q_1$ is a fibration,
$q_i\circ h=F_i, i=0,1.$ Passing to the corresponding derived categories
we get the functors $q_{i*}\circ p_{j*}: D(B)\to D(C)$ which are quasi-inverse
to $\alpha_*\circ\beta_*: D(C)\to D(B)$. This implies that $\theta_{F_0}=
\theta_{F_1}$.
\end{pf}

\subsection{Derived category of virtual modules}
\label{virtual}
\subsubsection{}
For a $\CO$-algebra $A$ we define the derived category
$D(\CO,A)$ to be the derived category of modules $DU(\CO,P)$ where
$P\to A$ is a cofibrant resolution. The category $D(\CO,A)$ is defined
uniquely in a way one could expect from an object of 2-category $\Cat$:
it is unique up to an equivalence which is unique up to a unique 
isomorphism.

Any map $f=(\alpha,\phi): (\CO,A)\to (\CO',A')$ of operad algebras
over a map $\alpha$ of operads defines a pair of adjoint functors
$$ \Left f^*:D(\CO,A)\rlarrows D(\CO',A'): \Right f_*$$
To construct these functors one has to choose cofibrant resolutions
$P\to A$ and $P'\to A'$ of the algebras; the $\CO'$-algebra $\phi^*(P)$
is cofibrant and therefore one can lift the composition
$$\phi^*(P)\to\phi^*(A)\overset{\alpha}{\to}A'$$
to a map $\phi^*(P)\to P'$. The construction is unique up to
a unique isomorphism.

We present below a more ``canonical'' construction
of $D(\CO,A)$ in terms of fibered categories. This approach follows
\cite{belu}, 2.4.

The correspondence $A\mapsto DU(\CO,A)$ together with the
functors $\Right f_*=f_*$ of~(\ref{adj-D}) as "inverse image functors"
define a fibered category $DU/\Alg(\CO)$. 

Let $A\in\Alg(\CO)$. Denote by $\COF/A$ the category of maps $P\to A$
with cofibrant $P$ and let $c_A:\COF/A\to\Alg(\CO)$ be the forgetful
functor defined by $c_A(P\to A)=P$.

\subsubsection{}
\begin{defn}{D(A)-fib} The derived category of virtual $(\CO,A)$-modules
is the fiber of $DU/\Alg(\CO)$ over $c_A$.
\end{defn}

In other words, an object of $D(\CO,A)$ consists of a collection
$X_a\in D(\CO,P)$ for each map $a:P\to A$ with cofibrant $P$, endowed with
compatible isomorphisms $\phi_f:X_a\to f_*(X_b)$ defined for any
presentation of $a$ as a composition
$$ P\overset{f}{\to}Q\overset{b}{\to}A.$$

\subsubsection{} Any object $\alpha:P\to A$ in $\COF/A$ defines an
obvious functor $q_{\alpha}:D(\CO,A)\to DU(\CO,P)$. Also the
functor $v_*:DU(\CO,A)\to D(\CO,A)$ is defined so that 
$q_{\alpha}\circ v_*$ and $\alpha_*$ are naturally isomorphic.

\begin{prop}{defs-coincide}
Let $\alpha:P\to A$ be a weak equivalence in $\COF/A$. Then the functor
$$ q_{\alpha}:D(\CO,A)\to DU(\CO,P)$$
is an equivalence. In part, the derived category of 
virtual $A$-modules  "is just" the derived category of modules over
the derived  enveloping algebra $\Left U(\CO,A)$.
\end{prop}
\begin{pf} We will omit the operad $\CO$ from the notations.
Let us construct a quasi-inverse functor $q^{\alpha}:DU(P)\to D(A)$.
For any $\beta:Q\to A$ in $\COF/A$ choose a map $f_{\beta}:Q\to P$.
This map is unique up to a homotopy $F:Q\to P^I$ for an appropriate
path diagram
$$P\overset{i}{\lra}P^I\rra P.$$
Moreover, the homotopy $F:Q\to P^I$ is itself unique up to a homotopy
as in~\ref{uniqueness-of-d-i}.
Now, for any $X\in DU(P)$ put $X(\beta)=(f_{\beta})_*(X)$. 
\Lem{uniqueness-of-d-i}(2) implies that the collection
of objects $\{X_{\beta}\}$ can be uniquely completed to an object
of $D(A)$.
\end{pf}

\Prop{defs-coincide} implies that the functor $v_*$ admits a left
adjoint functor $v^*:D(A)\to DU(A)$.

\subsubsection{} Let now  $f:A\to B$ be a morphism in $\Alg(\CO)$.
The functor $f_*:D(B)\to D(A)$ is induced by the obvious functor
$\COF/A\to\COF/B$.

\Prop{defs-coincide} implies that $f_*$ admits a left adjoint
functor $f^*:D(A)\to D(B)$.  Of course, it can be defined as
$$f^*=q^{\beta}\circ\Right g^*\circ q_{\alpha}$$
where $\alpha:P\to A$ and $\beta:Q\to B$ are cofibrant resolutions
of $A$ and $B$ respectively and a map $g:P\to Q$ satisfies the condition
$\beta\circ g=f\circ\alpha$.

\subsection{Varying the operad} 
Let now $\alpha:\CO\ra\CO'$ be a quasi-isomorphism of operads compatible with
$\Sigma$-splittings.
\subsubsection{}
\begin{thm}{comparison-5a}
There is an isomorphism of functors 
$$\Left U\lra\Left U\circ\Left\alpha^*.$$
\end{thm}
\begin{pf}
It suffices to prove that if $A$ is a cofibrant $\CO$-algebra, the 
composition
$$ U(\CO,A)\ra U(\CO,\alpha_*\alpha^*A)\ra U(\CO',\alpha^*A)$$
is a quasi-isomorphism.

The proof is similiar to that of~\ref{UA=UA-X}. The claim immediately
reduces to the case when $A$ is standart cofibrant and finitely generated.
Choose a free homogeneous base $\{x_i\}_{i\in I}$ for $A$; choose a
total order on $I$ so that $dx_i$ belongs to the subalgebra generated
by the elements $\{x_j\}_{j<i}$. Let $M$ be the set of multi-indices
ordered as in ~\ref{eq-opalg}. 

This defines filtrations on $U(\CO,A)$ and on $U(\CO',\alpha^*A)$ as 
in~(\ref{FU}).

The associated graded complexes take form
$$ \CO(|d|+1)\otimes_{\Sigma_d}k\text{ and }\CO'(|d|+1)\otimes_{\Sigma_d}k;$$
they are quasi-isomorphic since $\alpha$ is quasi-isomorphism preserving
$\Sigma$-splittings.
Theorem is proven.
\end{pf}

Putting together~\Cor{U(w-c)},~\Thm{asso-comparison} and~\Thm{comparison-5a}
we get immediately the following

\subsubsection{}
\begin{thm}{comparison-5}
Let $f=(\alpha,\phi): (\CO,A)\to (\CO',A')$ be a weak equivalence of
operad algebras.

Then the pair of derived functors
$$ \Left f^*:D(\CO,A)\rlarrows D(\CO',A'): \Right f_*$$
provides an equivalence of the derived categories.
\end{thm}

\section{Category of operads}
\label{section-catop}

\subsection{Introduction}
The category $\op(k)$ of operads has itself "algebraic" nature: an operad is
a collection of complexes endowed with a collection of operations
satisfying a collection of identities. This is why one can 
mimic the construction of Section~\ref{cmc} to define a CMC structure
on $\op(k)$.

The aim of this Section is to prove the following

\subsubsection{}
\begin{thm}{op-models} The category of operads $\op(k)$ in $C(k)$ admits
a structure of closed module category in which

--- $\alpha:\CO\ra\CO'$ is a weak equivalence if for all $n$ 
$\alpha_n$ is a quasi-isomorphism.

--- $\alpha:\CO\ra \CO'$ is a fibration if it is componentwise
surjective.
\end{thm}

The scheme of the proof is very close to the proof 
of~\Thm{cmc-structure}. In particular, a description of cofibrations in 
$\op(k)$ similar to that of~\Rem{retracts} will be given.

The proof of \Thm{op-models} is given in~\ref{free-operads}--
\ref{pf-op-models}. In \ref{standard-examples} we check that the standard
Lie and commutative operads $\CS$, $\CS_C$ from~\cite{hla} are cofibrant
operads in the sense of~\Thm{op-models}. Finally, in~\ref{more} we prove
that the derived category $D(\CO,A)$ can be ``calculated'' using a cofibrant
resolution of $\CO$ if $A$ is a flat $k$-complex.

\subsection{Free operads}
\label{free-operads}

The definitions below are close to~\cite{gk}, 1.1,2.1.

Let $\Col(k)$ be the category of {\em collections} of complexes in $C(k)$
numbered by nonnegative integers. As a category, this is a direct
product of $\Bbb N$ copies of $C(k)$. The obvious forgetful
functor $\#:\Op(k)\lra\Col(k)$ admits a left adjoint {\em free operad}
functor $F:\Col(k)\lra \Op(k)$ which can be described explicitly
using the language of trees.

\subsubsection{}
\begin{defn}{tree}(cf.~\cite{hla}, 4.1.3)
A tree is a finite directed graph with one initial (=having no ingoing edges)
vertex, such that any non-initial vertex has exactly one ingoing edge. 
\end{defn}

Terminal vertices of a tree are those having no outgoing edges; internal
vertices are those that are not terminal.

{\em Notations.} For a tree $T$ the set of its terminal (resp., internal)
vertices is denoted by $\ter(T)$ (resp., $\intl(T)$); $t(T)$ (resp., $i(T)$) 
is the number of terminal (resp., internal) vertices of $T$.  For any
$v\in\intl(T)$ the set of its outgoing vertices is denoted by $\out(v)$
and their number is $o(v)$.

{\em We choose once and forever a set of representatives of isomorphism classes
of trees; only these representatives will be called trees. 
}

For instance, for each $n$ we have a unique tree having one internal 
(=initial) vertex and $n$ terminal vertices. This is called $n$-corolla.

\subsubsection{}
\begin{defn}{ntree}
A $n$-tree consists of a pair $(T,e)$ where $T$ is a tree
and $e:\langle n\rangle\lra\ter(T)$ is an injective map.
\end{defn}

Denote by $\irr(T)$ the set of terminal vertices of $T$ which do not
belong to the image of $e$.

Denote by $\CT(n)$ the set of $n$-trees. 

The group $\Sigma_n$ acts on $\CT(n)$ on the right by the rule
$$ (T,e)\sigma=(T,e\sigma).$$

The collection $\CT=\{\CT(n)\}$ admits a structure of operad in the 
category $\Ens$. In fact, if $(T_0,e_0)\in\CT(n), (T_i,e_i)\in\CT(m_i)$
then the composition $T$ of the trees is defined by identifying the root
of $T_i$ with the terminal vertex $e(i)$ of $T$. The set $\ter(T)$
contains the disjoint union $\cup\ter(T_i)$ and the injective map
$e:\Sigma_m\lra \ter(T)$ for $m=\sum m_i$ is given by the formula
$$e(m_1+\ldots+m_{i-1}+j)=e_i(j)\in\ter(T_i)\subseteq\ter(T).$$

\subsubsection{}
\label{tree-free}
Here is the explicit construction of the free operad
functor. Let $V=\{V_i\}\in\Col(k)$. For any $T\in\CT(n)$ define a complex
$V_T$ by the formula
$$ V_T=\bigotimes_{v\in\intl(T)}V(o(v))\otimes V(0)^{\otimes\ter(T)-n}.$$
This should be interpreted as follows: each internal vertex $v$ of $T$
we mark with an element of $V(o(v))$; each non-numbered terminal vertex 
of $T$ we mark with an element of $V(0)$.

Note that for any $T\in\CT(n)$ and $\sigma\in\Sigma_n$ the complexes
$V_T$ and $V_{T\sigma}$ are tautologically isomorphic.

The free operad $F(V)$ generated by the collection $V$ is thus defined
by the formula
\begin{equation}
F(V)(n)=\bigoplus_{T\in\CT(n)} V_T.
\end{equation}
The $\Sigma_n$-action on $F(V)(n)$ is defined as follows. Let $x\in F_T,
\sigma\in\Sigma_n$. Then $x\sigma$ is "the same element as
$x$ but in $V_{T\sigma}$".

The operad multiplication is defined obviously by the multiplication in
the $\Ens$-operad $\CT$.

The map $V\lra\#F(V)$ of collections carries each $V(n)$ to the direct
summand $V_T$ of $F(V)(n)$ corresponding to the $n$-corolla
endowed with a(ny) bijective map $e:\langle n\rangle\ra\ter(T)$.

\subsection{Ideals; limits and colimits.}

Let $\CO\in\Op(k)$. An ideal $\CI$ in $\CO$ is a collection of $\Sigma_n$-
invariant subcomplexes $\{\CI(n)\subseteq\CO(n)\}$ which is stable
under the composition in an obvious way (if one of the factors belongs
to $\CI$ then the result belongs to $\CI$). A kernel of a map of operads
is always an ideal; if $\CI\subseteq\CO$ is an ideal then the quotient
operad $\CO/\CI$ is correctly defined.  If $X\subseteq\CO^{\#}$ is a
subcollection, the ideal $(X)$ is defined as the smallest ideal containing
$X$.

Limits in the category $\Op(k)$ exist and commute with the forgetful functor
$\#:\Op(k)\ra\Col(k)$. Colimits can be constructed using the free operad
construction: if $\alpha:I\ra\Op(k)$ is a functor, its colimit is the
quotient of the free operad generated by the collection 
$\dirlim\#\circ\alpha$ by an appropriately defined ideal.
Note that filtered colimits commute with $\#$.

\subsection{Adding a variable to kill a cycle}

Let $\CO$ be an operad in $C(k)$, $M\in C(k)$ and let $\alpha:M\ra\CO(n)$
be a map of complexes. The operad $\CO\langle M,n,\alpha\rangle$ is defined
as in~\ref{attaching-cmc}. If $M=k[d]$ and $\alpha:M\to\CO(n)$ takes the generator
of $M$ to a cycle $a\in\CO(n)^d$, the resulting operad is obtained by
``adding a variable to kill the cycle $a$''. It is denoted by 
$\CO\langle T; dT=a\rangle$.

One can immediately see that a map $\CO\to\CO\langle T; dT=a\rangle$ satisfies
the LLP with respect to any surjective quasi-isomorphism of operad. Also,
if the complex $M\in C(k)$ takes form $0\to k=k\to 0$, any map
$$\CO\to\CO\langle M,n,\alpha\rangle$$
satisfies the LLP with respect to any surjective map of operads.

Similarly to~\ref{standard:c,ac} one defines standard cofibrations and 
standard acyclic cofibrations as appropriate direct limits of the 
maps described.

\subsection{Extension of a homotopy to the free operad.}
\label{hfreeop}
Here we repeat the construction of~\ref{hfree}. In our case the construction
will be even easier since the operads are similar to associative algebras
and not to general operad algebras.

Let $\alpha:V\ra V$ be an endomorphism of a collection $V$ and
$h:V\ra V[-1]$ be a homotopy: $dh=\id_V-\alpha$. We wish to construct 
a homotopy $H:F(V)\ra F(V)[-1]$ between $\id_{F(V)}$ and $F(\alpha)$. 

For this we fix a total order on the set of terminal vertices of each 
corolla. This gives a lexicographic order on the set of all vertices
of any tree. The restriction of $H$ on $V_T$ is defined as
$$ H=\sum_{v\in\intl(T)\cup\irr(T)} H_v$$
where 
$$H_v=\bigotimes_{w\in\intl(T)\cup\irr(T)}\theta_w^v$$
with 
\begin{equation}
\theta_w^v=\begin{cases}
                \alpha\text{ if } w<v\\
                h\text{ if } w=v\\
                \id\text{ if } w>v
           \end{cases}
\end{equation}
One immediately checks that $dH=\id_{F(V)}-F(\alpha)$.
\subsubsection{}
\begin{lem}{H-invariance}
Let an ideal $\CI$ in the algebra $F(V)$ be generated by a set of elements
$\{x_i\}$. Then, if $H(x_i)\in\CI$ for all $i$, the ideal $\CI$ is
$H$-invariant.
\end{lem}
\begin{pf}
Straightforward calculation.
\end{pf}

\subsection{Proof of~\Thm{op-models}}
\label{pf-op-models}
The proof is close to that of~\Thm{cmc-structure} and~\Thm{opalg}.

Since $\Op(k)$ admits arbitrary limits and colimits, and the forgetful
functor $\#:\Op(k)\ra\Col(k)$ commutes with filtered colimits, we have
only to check that for any operad $\CO$ and contractible collection $X$
the natural map $\CO^{\#}\ra(\CO\coprod F(X))^{\#}$ is homotopy equivalence.

We proceed as in the proof of~\Thm{opalg}. Put 
$V=\CO^{\#}\oplus X\in\Col(k)$. If $I$ is the kernel of the natural
projection $F(\CO^{\#})\lra\CO$ and $J$ is the ideal in $F(V)$ generated
by $I$, then $\CO\coprod F(X)$ is isomorphic to $F(V)/J$.

Let $\alpha:V\ra V$ be the composition $V=\CO^{\#}\oplus X\ra\CO^{\#}\ra V$
and let $h:V\ra V[-1]$ be the homotopy, $dh=\id_V-\alpha$, vanishing on 
$\CO^{\#}$. According~\ref{hfreeop} a homotopy $H:F(V)\ra F(V)[-1]$ is
defined and by~\Lem{H-invariance} the ideal $J$ is $H$-invariant.
Then $H$ induces a homotopy on $F(V)/J=\CO\coprod F(X)$ and this proves the
theorem.

\subsection{Standard examples.}
\label{standard-examples}

Suppose that $k$ contains $\Bbb Q$. Recall that dg Lie algebras 
(resp., commutative dg agebras) are precisely algebras over an
appropriate operad $\Lie$ (resp., over $\Com$). Their strong homotopy
counterparts are correspondingly the algebras over the  "standard" operads 
$\CS$ and $\CS_C$ see~\cite{hla}, 4.1 and 4.4.

Let us show that the standard operads,
$\CS$ (the standard Lie operad) and $\CS_C$ (the standard commutative
operad) are cofibrant in $\Op(k)$.

These operads are constructed by consecutive "attaching a variable to 
kill a cycle" ---
as it is explained in {\em loc. cit.}, 4.1.1., which differs a little
from our construction. However, if $k\supseteq\Bbb Q$, this operation 
also gives rise to a cofibration as shows the~\Lem{attaching} below.

Let $\CO$ be an operad and let $z\in\CO(n)$ be a cycle. Let $G\subseteq
\Sigma_n$ and a character $\chi:G\ra k^*$ satisfy the condition
$zg=\chi(g)z$ for all $g\in G$.  Then the "attaching of a variable"
is defined (well, is not defined) in {\em loc. cit.} to be the map
$$\CO\lra\CO'=\CO\langle e;de=x\rangle/(eg-\chi(g)e).$$
\subsubsection{}
\begin{lem}{attaching}
In the notations above the map $\CO\ra\CO'$ is a cofibration.
\end{lem}
\begin{pf}
The projection $\CO\langle e;de=x\rangle\lra\CO'$ is split by the map
$\CO'\ra\CO\langle e;de=x\rangle$ which sends the element $e$ to
$\frac{1}{|G|}\sum_{g\in G}\chi(g^{-1})eg.$

Thus, the map $\CO\ra\CO'$ is a retract of a standard cofibration and the
Lemma is proven.
\end{pf}

\subsection{More on the derived category}
\label{more}

In \ref{defs-coincide} we saw that the derived category $D(\CO,A)$ of
virtual $A$-modules can be calculated using a cofibrant resolution
of $A$ in $\Alg(\CO)$. Now we will show that one can take a cofibrant
resolution of the operad $\CO$ instead. 

\subsubsection{}
\begin{lem}{U:c(o)}Let $\CO$ be an operad over $k$, $\alpha:M\to\CO(n)$
be a map of complexes. Put $\CO'=\CO\langle M,n,\alpha\rangle$.
Let $A$ be a $\CO'$-algebra and let $U=U(\CO,A)$, $U'=U(\CO',A)$.
Then one has a natural isomorphism
$$ U'\iso U\langle M\otimes A^{\otimes n-1},\tilde{\alpha}\rangle$$
where $\tilde{\alpha}$ is the composition
$$ M\otimes A^{\otimes n-1}\overset{\alpha\otimes\id}{\lra}
\CO(n)\otimes A^{\otimes n-1}\to U(\CO,A).$$   
\end{lem}
\begin{pf}An $\CO'$-algebra structure on a $\CO$-algebra $A$ is given 
by a map $f:M\otimes A^{\otimes n}\to M$ such that $d(f)$ is equal
to the composition
$$ M\otimes A^{\otimes n}\overset{\alpha\otimes\id}{\lra}\CO(n)\otimes 
A^{\otimes n}\to A.$$

A structure of $(\CO',A)$-module on a $(\CO,A)$-module $X$ is given by 
a map $m:M\otimes A^{\otimes n-1}\otimes X\to X$ such that $d(m)$ is 
equal to the composition
$$ M\otimes A^{\otimes n-1}\otimes X\overset{\alpha\otimes\id\otimes\id}
{\lra}\CO(n)\otimes A^{\otimes n-1}\otimes X\to X.$$  
This proves the claim. 
\end{pf}

The following \Lem{asso+flat} will be used in~\ref{U:co(quis)}.

\subsubsection{}
\begin{lem}{asso+flat}
Let $A,B\in\dga(k)$, $M,N\in C(k)$ and let a commutative diagram

\begin{center}
  $\begin{array}{ccc}
     M&\overset{\alpha}{\lra}& A\\
     \Big\downarrow\vcenter{\rlap{$\scriptstyle g$}}&&
     \Big\downarrow\vcenter{\rlap{$\scriptstyle f$}} \\
     N&\overset{\beta}{\lra}& B
  \end{array}   $
\end{center}

be given so that $f:A\to B$ is a weak equivalence and $g:M\to N$ is
a quasi-isomorphism. If $A$, $B$, $M$, $N$ are flat $k$-complexes then the
induced map
$$ A\langle M,\alpha\rangle\to B\langle N,\beta\rangle$$
is a weak equivalence. Moreover, the algebras
$ A\langle M,\alpha\rangle, B\langle N,\beta\rangle$ are also
flat over $k$.
\end{lem}
\begin{pf}
The associative algebra $A\langle M,\alpha\rangle$ admits a natural
filtration $\{F_n\}$ defined by
$$ F_n=\sum_{k=0}^n (A\otimes M[1])^{\otimes k}\otimes A.$$
The associated graded pieces are
$$ \gr_n= (A\otimes M[1])^{\otimes n}\otimes A.$$
Now it is clear that the map in question induces isomorphism
of the associated graded pieces and therefore is itself a 
quasi-isomorphism. The complex $F_n$ can be obtained  as the cone
of a map $\gr_n\to F_{n-1}$ induced by $\alpha$. 

Since flatness is closed under taking cones
and filtered colimits (see~\Lem{cofibrant->flat}), the new algebras 
$ A\langle M,\alpha\rangle, B\langle N,\beta\rangle$ are
flat $k$-complexes.
\end{pf}

\subsubsection{}
\begin{cor}{U:co(quis)}Let $\CO$ be a cofibrant operad and let 
$f:A\to A'$ be a quasi-isomorphism of $\CO$-algebras. If $A$ and $A'$
are flat as complexes over $k$ then the natural map  $U(f)$ is 
quasi-isomorphism.
\end{cor} 
\begin{pf} One can suppose $\CO$ to be standard cofibrant.
Then the claim follows immediately from~\ref{U:c(o)} and~\ref{asso+flat}.
\end{pf}

\subsubsection{}
\begin{prop}{D=DU}
Let $A$ be an algebra over a cofibrant operad 
$\CO\in\op(k)$.
Then, if $A$ is flat as a $k$-complex, there exists a natural equivalence
$$D(\CO,A)\to DU(\CO,A).$$
\end{prop}
\begin{pf}The claim will immediately follow from~\ref{U:co(quis)}
once we check that cofibrant algebras over cofibrant operads are flat.

This immediately reduces to the case of a finitely generated
standard cofibrant algebra. This one admits a filtration $\{\CF_n\}$
as in~\ref{fg-cof}. Thus everything is reduced to checking that
for any cofibrant operad $\CO$ the complex $\CO(n)\otimes_{\Sigma_n}k$
is flat. This follows from the tree description of a free operad
--- see~\ref{tree-free}.

\end{pf}

\subsubsection{}
\begin{exa}{Lie}(see the notations of~\ref{standard-examples})
Let $A$ be a flat dg Lie algebra over 
$k\supseteq\Bbb{Q}$. Then the derived category $DU(\Lie,A)$ of $A$-modules 
is equivalent to the derived category $DU(\CS,A)$ of modules over
$A$ considered as a strong homotopy Lie algebra. In fact,
the category $DU(\CS,A)$ is equivalent to $D(\CS,A)$ 
by~\Prop{D=DU} and \ref{standard-examples}.
The category $DU(\Lie,A)$ is equivalent to $D(\Lie,A)$ by the
PBW theorem for dg Lie algebras --- the latter implies that enveloping
algebras of quasi-isomorphic flat Lie algebras are quasi-isomorphic.
Finally, the categories
$D(\CS,A)$ and $D(\Lie,A)$ are naturally equivalent by~\Thm{comparison-5}.
\end{exa}

\section{Cotangent complex; cohomology of operad algeras.}
\label{section-cot}

\subsection{Introduction}
Let $\CO$ be a $\Sigma$-split operad in $C(k)$. In this Section we construct
for an algebra map $B\to A$ in $\Alg(\CO)$ its cotangent complex 
$L_{A/B}\in D(\CO,A)$ belonging to the derived category of virtual
$A$-modules. 

In the next Section we define the tangent complex $T_A$ which is, as usual,
dual to the cotangent complex. The tangent complex
admits a unique (in the homotopy category) structure of dg Lie algebra.
This Lie algebra must play a crucial role in the deformation theory of 
operad algebras.

\subsection{Derivations} 
\subsubsection{}
\begin{defn}{der}
Let $\CO$ be an operad in $C(k)$, $\alpha: B\to A$ be a map in $\Alg(\CO)$
and let $M$ be a $A$-module. A map $f:A\to M$ of complexes (not necessarily
commuting with the differentials) is called {\em $\CO$-derivation over $B$} 
if it vanishes on the image of $B$ in $A$ and for each $n>0$ the 
following diagram
\begin{center}
\begin{picture}(10,4)
  \put(0,0){\makebox(5,1){$\underset{a+b=n-1}{\sum}\CO(n)\otimes 
       A^{\otimes a}\otimes M\otimes A^{\otimes b}$}}
  \put(9,3){\makebox(1,1){$A$}}
  \put(0,3){\makebox(5,1){${\cal O}(n)\otimes A^{\otimes n}$}}
  \put(9,0){\makebox(1,1){$M$}}
  \put(2.5,3){\vector(0,-1){2}}
  \put(9.5,3){\vector(0,-1){2}}
  \put(6.3,0.5){\vector(1,0){2}}
  \put(4.5,3.5){\vector(1,0){4}}
  \put(9.5,1){\makebox(.5,2){$\scriptstyle f$}}
  \put(3,1){\makebox(3.5,2){$\scriptstyle \underset{a+b=n-1}{\sum}
        1\otimes 1^a\otimes f\otimes 1^b$}}
\end{picture}
\end{center}
is commutative.
\end{defn}
\subsubsection{}
The $\CO$-derivations from $A$ to $M$ over $B$ form a complex 
$\Der^{\CO}_{B}(A,M)$. When $B=F(0)=\CO(0)$ is the initial object in 
$\Alg(\CO)$ we will omit the subscript $B$ from the notation.
We will also omit the superscript $\CO$ when it does not make a confusion.

The complex of derivations $\Der^{\CO}_{B}(A,M)$
is a subcomplex of $\CHom_k(A,M)$. This defines a functor
$$ \Der^{\CO}_B(A,\_): \Mod(A)\to C(k)$$
which is representable in the following sense.

\begin{prop}{omega}
There exists a (unique up to a unique isomorphism) $A$-module $\Omega_{A/B}$
(called the module of relative differentials) together with a derivation 
$\partial:A\to\Omega_{A/B}$ inducing the natural isomorphism of complexes
$$ \CHom_A(\Omega_{A/B},M)\overset{\sim}{\lra}\Der^{\CO}_B(A,M).$$
\end{prop}
\begin{pf} 1. Consider firstly the absolute case $B=F(0)$.
The functor $M\mapsto\CHom_k(A,M)$ is obviously represented
by the free $A$-module $U(\CO,A)\otimes A$. Thus, the functor $\Der(A,\_)$
is represented by the quotient of $U(\CO,A)\otimes A$ modulo the relations
which guarantee the commutativity of the diagrams

\begin{center}
\begin{picture}(10,4)
  \put(0,0){\makebox(5,1){$\underset{a+b=n-1}{\sum}\CO(n)\otimes 
     A^{\otimes a}\otimes \left( U(\CO,A)\otimes A \right)     
     \otimes A^{\otimes b}$}}
  \put(9,3){\makebox(1,1){$A$}}
  \put(0,3){\makebox(5,1){${\cal O}(n)\otimes A^{\otimes n}$}}
  \put(9,0){\makebox(1,1){$U(\CO,A)\otimes A$}}
  \put(2.5,3){\vector(0,-1){2}}
  \put(9.5,3){\vector(0,-1){2}}
  \put(7,0.5){\vector(1,0){1.1}}
  \put(4.5,3.5){\vector(1,0){4}}
  \put(9.5,1){\makebox(.5,2){$\scriptstyle f$}}
  \put(3,1){\makebox(3.5,2){$\scriptstyle \sum\limits_{a+b=n-1}^{}
           1\otimes 1^a\otimes f\otimes 1^b$}}
\end{picture}
\end{center}
Here $f:A\ra U(\CO,A)\otimes A$ is given by $f(a)=1\otimes a$.

2. In general one sees immediately that the complex
$\Omega_A/\partial\circ\alpha(B)$ represents the functor $\Der^{\CO}_B(A,\_)$.
\end{pf}   

\subsubsection{} Let $C\overset{\alpha}{\to}B\overset{f}{\to} A$ be a pair
of maps in $\Alg(\CO)$. Any derivation
$d:A\to M$ over $C$ defines a derivation $d\circ f:B\to f_*(M)$ over $C$. 
This defines a canonical map 
$\Omega^f:f^*\Omega_{B/C}\to\Omega_{A/C}$. 

On the other hand, any derivation $d:A\to M$ over $B$ can be considered
as a derivation over $C$. This defines a canonical map
$\Omega_{\alpha}:\Omega_{A/C}\to\Omega_{A/B}$.

\Prop{more-corr-cot} below claiming that sometimes $\Omega^f$ and 
$\Omega_{\alpha}$ are quasi-isomorphisms, allows one to define correctly
the cotangent complex $L_{A/B}\in D(\CO,A)$ as a (sort of) left
derived functor of the functor $\Omega_{A/B}$ --- see~\ref{cot}.

\subsubsection{}
\label{?varop}
More generally, given a map $\alpha:\CO\to\CO'$ of operads and a couple 
of maps $f:B\to A\in\Mor\Alg(\CO)$, $f':B'\to A'\in\Mor\Alg(\CO')$,  
one defines $u:f\to f'$ over $\alpha$ to be the pair of maps
$\phi:\alpha^*B\to B',\psi:\alpha^*A\to A'$
making the diagram
\begin{center}
$\begin{array}{ccc}
        \alpha^*B&\overset{\phi}{\lra}& B'\\
        \Big\downarrow\vcenter{\rlap{$\scriptstyle f$}}&&
        \Big\downarrow\vcenter{\rlap{$\scriptstyle f'$}} \\
        \alpha^*A&\overset{\psi}{\lra}& A'
  \end{array}   $
\end{center}
commutative.

For any $A'$-module $M$ one has a canonical ``restriction map''
$$ \Der^{\CO'}_{B'}(A',M)\to\Der^{\CO}_B(A,f_*(M)) $$
which gives a canonical map
$ \Omega^u: u^*(\Omega_{A/B})\to\Omega_{A'/B'}.$
Of course, $\Omega^u$ can be presented as a composition of three maps ---
the map $\Omega^{\alpha}:\Omega_{A/B}\to\Omega_{\alpha^*A/\alpha^*B}$
responsible for the change of operads 
($\phi=\id_{\alpha^*B},\psi=\id_{\alpha^*A}$); the map 
$\Omega^{\psi}:\Omega_{\alpha^*A/\alpha^*B}\to\Omega_{A'/\alpha^*B}$
and the map
$\Omega_{\phi}:\Omega_{A'/\alpha^*B}\to\Omega_{A'/B'}$.

\subsection{Module of differentials of a cofibration}

\subsubsection{} Let  $A=F(X)$ be the 
free $\CO$-algebra generated by a complex $X\in C(k)$. Then there is a 
natural isomorphism
$$ \Omega_A=U(\CO,A)\otimes X.$$
We wish to describe the module $\Omega_{A/B}$ when $\alpha:B\to A$ is a
standard cofibration.

\subsubsection{}
Let $\alpha:B\to A$ be a map of $\CO$-algebras, $M\in C(k)$.
Let $f:M\to A^{\#}$ be a map in $C(k)$ and let 
$A'=A\langle M,f\rangle$ be defined as in~\ref{attaching-cmc}.
Put $U=U(\CO,A),\ U'=U(\CO,A')$. The map $\partial\circ\alpha:M\to\Omega_{A/B}$
defines $\alpha':U'\otimes M\to U'\otimes_U\Omega_{A/B}$.

\begin{lem}{adj-omega} The $A'$-module $\Omega_{A'/B}$
is naturally isomorphic to the cone of $\alpha'$.
\end{lem}
\begin{pf} Any $B$-derivation from $A'$ to a $A'$-module is uniquely defined by
its restrictions to $A$ and to $M$ satisfying the obvious compatibility
condition including the map $f$. Thus $\Omega_{A'/B}$ and $\cone(\alpha')$ 
represent the same functor. 
\end{pf}

\subsubsection{}
\begin{cor}{omega-is-cof}
Let $\alpha:B\to A$ be is a cofibration. Then $\Omega_{A/B}$ is a cofibrant
$A$-module (in the sense of~\ref{dha1-models}).
\end{cor}
\begin{pf}
If $\alpha$ is a standard cofibration, \Lem{adj-omega} immediately implies that
$\Omega_A$ is semi-free, i.e., standard cofibrant in $\Mod(A)$.

In the general case, let $i:A\rlarrows C: p$ represent $\alpha$ as a retract
of a standard cofibration $i\circ\alpha:B\to C$. Then the maps
$$\Omega_{A/B}=p^*i^*\Omega_{A/B}\to p^*\Omega_{C/B}\to \Omega_{A/B}$$
define $\Omega_{A/B}$ as a retract of $p^*\Omega_{C/B}$. Since the inverse
image functor preserves cofibrations, the claim follows.
\end{pf}

\subsubsection{}
\begin{prop}{compos+cof}
Let $C\overset{\alpha}{\to}B\overset{f}{\to} A$ be a pair
of maps in $\Alg(\CO)$ so that $f$ is a cofibration. Then the
sequence
\begin{equation}
0\lra f^*\Omega_{B/C}\overset{\Omega^f}{\lra}\Omega_{A/C}
\overset{\Omega_{\alpha}}{\lra}\Omega_{A/B}\lra 0
\label{ses-omega}
\end{equation}
is exact.
\end{prop}
\begin{pf}
If $f$ is a standard cofibration one proves the claim by induction
using~\Lem{adj-omega}.

To prove the general case, let $A'\overset{q}{\to}A\overset{j}{\to}A'$
satisfy $q\circ j=\id_A$ so that $f'=j\circ f$ is a standard cofibration.
If $\fs$ denotes the sequence~(\ref{ses-omega}) and $\fs'$ the
same sequence constructed for $f'$ instead of $f$, one immediately
obtains that $j^*(\fs)$ is a retract of $\fs'$ and is, therefore, exact.
Since the $A$-module $\Omega_{A/B}$ is cofibrant, $j^*(\fs)$ splits
and so it remains exact after application of $q^*$. Therefore, 
$\fs=q^*j^*(\fs)$ is exact.
\end{pf}

\subsubsection{}
\begin{cor}{corr-cot-cof}
Let $C\overset{\alpha}{\to}B\overset{f}{\to} A$ be a pair
of maps in $\Alg(\CO)$ where $f$ is an acyclic cofibration. Then the
natural map $\Omega^f$ is a quasi-isomorphism. 
\end{cor}
\begin{pf}
If $f$ is a standard acyclic cofibration, the first claim follows immediately
from~\Lem{adj-omega}. If $f$ is any acyclic cofibration, it is a retract of
a standard acyclic cofibration and the claim follows from the axiom~(CM 3).
\end{pf}

\subsubsection{}
\begin{prop}{corr-cot}
Let $C\overset{\alpha}{\to}B\overset{f}{\to} A$ be a pair
of maps in $\Alg(\CO)$. If $f$ is a quasi-isomorphism and $\alpha,f\circ\alpha$ 
are cofibrations then the map $\Omega^f$ is a weak equivalence.
\end{prop}
\begin{pf} Since $\alpha$ and $f\circ\alpha$ are cofibrations, 
there exists a map $g:A\to B$ over $C$ homotopy inverse to $f$.
Let $A\overset{i}{\to} A^I\overset{p_0,p_1}{\rra}A$ be a path 
diagram for $A$ such that $i$ is an acyclic cofibration and
the homotopy between $\id_A$ and $f\circ g$ is given by a map 
$h:A\to A^I$.

Since $i$ is an acyclic cofibration the map $i^*\Omega_{A/C}\to 
\Omega_{A^I/C}$
is a quasi-isomorphism of cofibrant $A^I$-modules; thus for $i=0,1$ the
maps $p^*_i\Omega_{A^I}\to \Omega_A$ are quasi-isomorphisms as well. 
The map
$\Omega^h:h^*\Omega_{A/C}\to\Omega_{A^I/C}$ becomes therefore 
quasi-isomorphism of cofibrant modules after application of $p_1^*$. 
Since $p_1$ is a weak equivalence of cofibrant algebras, \Cor{U(w-c)}
implies that $\Omega^h$ itself is a quasi-isomorphism. 
Then $p_2^*(\Omega_h)$ is also quasi-isomorphism. Thus, the map
$\Omega^{fg}:(fg)^*\Omega_{A/C}\to\Omega_{A/C}$ is a quasi-isomorphism.

In the same way, the map $\Omega^{gf}:(gf)^*\Omega_{B/C}\to\Omega_{B/C}$
is a quasi-isomorphism. This immediately implies that $\Omega^f$ is 
isomorphism.
\end{pf}

\subsubsection{}
\begin{prop}{more-corr-cot}
Let 
\begin{center}
$\begin{array}{ccc}
        B&\overset{\alpha}{\lra}& B'\\
        \Big\downarrow\vcenter{\rlap{$\scriptstyle f$}}&&
        \Big\downarrow\vcenter{\rlap{$\scriptstyle f'$}} \\
        A&\overset{\beta}{\lra}& A'
  \end{array}   $
\end{center}
be a commutative diagram in $\Alg(\CO)$ so that $\alpha,\beta$ are weak 
equivalences, $f,f'$ are cofibrations and $B,B'$ are cofibrant.
Then the composition
$$ \beta^*\Omega_{A/B}\overset{\Omega^{\beta}}{\lra}\Omega_{A'/B}
\overset{\Omega_{\alpha}}{\lra}\Omega_{A'/B'}$$
is a weak equivalence.
\end{prop}
\begin{pf} Applying \Prop{compos+cof} thrice we get the following
commutative diagram
\begin{center}
$\begin{array}{ccccccccc}
0&\to&\beta^*f^*\Omega_B&\to&\beta^*\Omega_A&\to&\beta^*\Omega_{A/B}&\to&0\\
&&        \Big\downarrow && \Big\downarrow && \Big\downarrow \\
0&\to&\beta^*f^*\Omega_B&\to&\Omega_{A'}&\to&\Omega_{A'/B}&\to&0\\ 
&&        \Big\downarrow && \Big\downarrow && \Big\downarrow \\
0&\to&g^*\Omega_{B'}&\to&\Omega_{A'}&\to&\Omega_{A'/B'}&\to&0\\
  \end{array}   $
\end{center}
Its rows are split exact sequences and the first two columns are 
quasi-isomorphisms (or even isomorphisms). Therefore the last column
we need consists of quasi-isomorphisms as well.
\end{pf}

\Lem{adj-omega} implies easily the following
\subsubsection{}
\begin{cor}{varop-omega}
Let $\alpha:\CO\to\CO'$ be a map of operads, $f\in\Mor\Alg(\CO)$,
$f'\in\Mor\Alg(\CO')$, $u:f\to f'$ be as in~\ref{?varop}.
If $\alpha,\phi,\psi$ are weak equivalences, $B,B'$ are cofibrant
and $f,f'$ are cofibrations then $\Omega^u$ is a quasi-isomorphism.
\end{cor}
\begin{pf} By~\ref{more-corr-cot} one can suppose that 
$B'=\alpha^*(B), A'=\alpha^*(A)$, so that we have to check that
$\Omega^{\alpha}$ is a weak equivalence.  Actually we shall prove
that this is isomorphism.

If $f:B\to A$ is a standard cofibration, the result immediately follows 
from~\Lem{adj-omega} by induction. Otherwise, if $f$ is a retract
of a standard cofibration $g:B\to C$, then $\Omega_{A/B}$
is a retract of $\pi^*\Omega_{C/B}$ where $\pi:C\to A$ is the
corresponding projection. This proves the claim.
\end{pf}

\subsection{Cotangent complex}

\subsubsection{Model structure on $\MOR(\CC)$}
\label{mod-mor}
Let $\CC$ be a closed model category. Then the category $\MOR(\CC)$
of arrows in $\CC$ admits a closed model category as follows.

A map from $g:Q\to P$ to $f:B\to A$ given by a pair of maps
$\beta:Q\to B$ and $\alpha:P\to A$ satisfying the condition
$f\circ\beta=\alpha\circ g$ is:

--- a weak equivalence (resp., a fibration) if both $\alpha$ and $\beta$ 
are weak equivalences (resp., fibrations);

--- a cofibration if $\beta$ and also the natural map $B\coprod^QP\to A$
are cofibrations.

We will use this model structure for the case $\CC=\Alg(\CO)$.
Note that, in particular, a map $g:Q\to P$ is cofibrant if $Q$ is a cofibrant
object and $g$ is a cofibration. The following lemma is the result of
the existence of the above defined model structure on $\MOR(\Alg(\CO))$.

\subsubsection{}
\begin{lem}{cof-res-arrow-contra}
1. For any map $f:B\to A$ in $\Alg(\CO)$ there exists a 
cofibrant resolution.

2. For any pair $g:Q\to P$, $g':Q'\to P'$ of cofibrant resolutions
of $f$ there exists a map, unique up to homotopy, from $g'$ to $g$.
\end{lem}

In the definition of cotangent complex below we use the notations of
\ref{defs-coincide}.
\subsubsection{}
\begin{defn}{cot}Let $f:B\to A$ be an algebra morphism in $\Alg(\CO)$.
The cotangent complex of $f$, $L_{A/B}$, is the object of $D(\CO,A)$
defined by the formula
$$ L_{A/B}=q^{\alpha}(\Omega_{P/Q})$$
where $g:Q\to P$ together with $\beta:Q\to B,\ \alpha:P\to A$
define a cofibrant resolution of $f$.
\end{defn}

\subsubsection{}
\begin{prop}{cot-correct}
The cotangent complex $L_{A/B}$ is defined uniquely up to a
unique isomorphism. 
\end{prop}
\begin{pf}
This immediately follows from \Lem{cof-res-arrow-contra} and 
\Prop{defs-coincide}.
\end{pf}

\subsubsection{Functoriality}

Let $\alpha:\CO\to\CO'$ be a map of operads, $f:B\to A$ and $f':B'\to A'$
be maps in $\Alg(\CO)$ and in $\Alg(\CO')$ respectively, and let
$u:f\to f'$ be a map over $\alpha$ as in~\ref{?varop}.

If $g:Q\to S$ is a cofibrant resoution of $f$, and $g':Q\to P'$ is a cofibrant
resolution for $f'$, there is, according to~\ref{mod-mor}, a map
$v:\alpha^*(g)\to g'$ aking the corresponding diagram commutative; moreover,
this map is unique up to homotopy. This defines a map
$$ L^u: \Left u^*(L_{A/B})\to L_{A'/B'}.$$

\Cor{varop-omega} immediately gives the following

\begin{prop}{varop-cot}
Let $\alpha:\CO\to\CO'$ be a map of operads, $f\in\Mor\Alg(\CO)$,
$f'\in\Mor\Alg(\CO')$, $u:f\to f'$ be as in~\ref{?varop}.
If $\alpha,\phi,\psi$ are weak equivalences then 
$L^u:\Left u^*(L_{A/B})\to L_{A'/B'}$ is an isomorphism in $D(\CO',A')$. 
\end{prop}

\subsection{Cohomology}

Let $A$ be a $\CO$-algebra, $M\in D(\CO,A)$ be a virtual $A$-module. 
The (absolute) cohomology of $A$ with coefficients in $M$ are defined to be
$$ H(A,M)=\Right\Hom_A(L_A,M)\in D(k).$$

\Prop{varop-cot} immediately implies the following comparison theorem

\subsubsection{}
\begin{thm}{comparison:opalg}
Let $f=(\alpha,\phi):(\CO',A')\to(\CO,A)$ be a weak equivalence of operad
algebras. Let $M\in\Mod(\CO,A), M\in\Mod(\CO',A')$ and let $g:M'\to f_*(M)$ 
be a quasi-isomorphism of $A'$-modules. Then the induced map
$$ H(A',M')\to H(A,M)$$
is an isomorphism in $D(k)$.
\end{thm}

\section{Tangent Lie algebra}
\label{section-tan}

Let $\CO$ be a $\Sigma$-split operad. For a cofibrant $\CO$-algebra
$A$ its tangent Lie algebra is defined to be $T_A=\Der^{\CO}(A,A)$.

The aim of this Section is to extend this correspondence to a functor  
from the category $\hoa(\CO)^{\is}$ of homotopy $\CO$-algebras and 
isomorphisms to the category $\hol(k)^{\is}$ of homotopy Lie algebras and
isomorphisms.

\subsection{} For any map $\alpha:A\to B$ let $\Der_{\alpha}(A,B)$
be the complex of derivations from $A$ to $B$ considered as a $A$-module
via $\alpha$. One has a pair of maps
$$ T_A\overset{\alpha_*}{\to}\Der_{\alpha}(A,B)
\overset{\alpha^*}{\leftarrow} T_B.$$

\begin{lem}{}
Let $A$ and $B$ be cofibrant $\CO$-algebras and let $\alpha$ be a weak
equivalence. Then $\alpha^*$ and $\alpha_*$ are quasi-isomorphisms.
\end{lem}
\begin{pf}
Recall that 
$$T_A=\CHom_A(\Omega_A,A),\ T_B=\CHom_B(\Omega_B,B),
\Der_{\alpha}(A,B)=\CHom_A(\Omega_A,B).$$
Since $A$ and $B$ are cofibrant, $\Omega_A$ and $\Omega_B$
are cofibrant modules over $A$ and $B$ respectively. The map $\alpha_*$
is a weak equivalence since $\Omega_A$ is cofibrant and $\alpha$ is
a quasi-isomorphism. The map $\alpha^*$ is a weak equivalence since
$\Omega^{\alpha}:\alpha^*(\Omega_A)\to\Omega_B$ is a weak equivalence
of cofibrant $B$-modules by~\Prop{corr-cot}.
\end{pf}

\subsection{Acyclic fibrations}

Let $\alpha:A\to B$ be an acyclic fibration (= surjective 
quasi-isomorphism). Put $I=\Ker \alpha$. Define
$$ T_{\alpha}=\{\delta\in T_A|\delta(I)\subseteq I\}.$$
Then $T_{\alpha}$ is a dg Lie subalgebra of $T_A$ and a natural Lie
algebra map $\pi_{\alpha}:T_{\alpha}\to T_B$ is defined.
Denote by $\iota_{\alpha}:T_{\alpha}\to T_A$ the natural inclusion.

\subsubsection{}
\begin{prop}{T(af)}The map $\pi_{\alpha}$ is a surjective quasi-isomorphism
while $\iota_{\alpha}$ is an injective quasi-isomorphism.
\end{prop}
\begin{pf}
Step 1. Let us check that $\pi_{\alpha}$ is surjective. Suppose first 
of all that $A$ is standard cofibrant i.e. is obtained from nothing
by a successive joining of free variables. Derivation on $A$ is uniquely
defined by its values on the free generators therefore any derivation
on $B$ can be lifted to a derivation on $A$ and it will belong
automatically to $T_{\alpha}$. For a general cofibrant $A$  let $C$
be a standard cofibrant algebra so that $A$ is a retract of $C$.
This means that there are maps $i:A\to C$ and $p:C\to A$ so that 
$p\circ i=\id_A$. Put $J=\Ker(\alpha\circ p)$. Then if $\delta\in T_B$
and if $\tilde{\delta}\in T_C$ lifts $\delta$ then the composition
$p\circ\delta\circ i$ is a derivation of $A$ which lifts $\delta$.

Step 2. One has 
$$\Ker\pi_{\alpha}=\{\delta\in T_A|\delta(A)\subseteq I\}=\Der(A,I).$$
Since $I$ is contractible, $\Ker\pi_{\alpha}$ is contractible.

Taking into account Steps 1 and 2 we deduce that $\pi_{\alpha}$ is
a surjective quasi-isomorphism.

Step 3. The diagram
\begin{center}
   \begin{picture}(4,3)
	\put(0,0){\makebox(1,1){$T_{\alpha}$}}
	\put(0,2){\makebox(1,1){$T_A$}}
	\put(2,0){\makebox(2,1){$T_B$}}
	\put(2,2){\makebox(2,1){$\Der_{\alpha}(A,B)$}}
	\put(0.5,1){\vector(0,1){1}}
	\put(3,1){\vector(0,1){1}}
	\put(1,0.5){\vector(1,0){1.4}}
	\put(1,2.5){\vector(1,0){0.7}}

	\put(1,2.5){\makebox(0.7,0.5){$\scriptstyle\alpha_*$}}
	\put(1,.5){\makebox(1.4,0.5){$\scriptstyle\pi_{\alpha}$}}
	\put(0.5,1){\makebox(0.7,1){$\scriptstyle\iota_{\alpha}$}}
	\put(3,1){\makebox(0.7,1){$\scriptstyle\alpha^*$}}
   \end{picture}
\end{center}
is commutative. Since the maps $\alpha^*,\alpha_*,\pi_{\alpha}$ are
quasi-isomorphism, $\iota_{\alpha}$ is also quasi-isomorphism.
Proposition is proven.
\end{pf}

\Prop{T(af)} allows one to define the map $T(\alpha): T_A\to T_B$
in the homotopy category $\hol(k)$ as 
$$T(\alpha)=\pi_{\alpha}\circ\iota_{\alpha}^{-1}.$$

\subsubsection{}
\begin{lem}{comp-af} Let $A\overset{\alpha}{\to}B\overset{\beta}{\to}C$
be a pair of acyclic fibrations. Then one has
$$ T(\beta\circ\alpha)=T(\beta)\circ T(\alpha)$$
in $\hol(k)$.
\end{lem}
\begin{pf}
Put $T=T_{\alpha}\times_{T_B}T_{\beta}.$ The map $T\to T_{\beta}$
is a surjective quasi-isomorphism since it is obtained by a base change
from $\pi_{\alpha}$. The Lie algebra $T$ identifies with
$$ \{\delta\in T_A|\delta(\Ker\alpha)\subseteq\Ker\alpha\text{ and }
\delta(\Ker(\beta\circ\alpha))\subseteq\Ker(\beta\circ\alpha)\}.$$
Therefore, $T$ is a subalgebra of $T_{\beta\circ\alpha}$ and all the
maps involved are quasi-isomorphisms. This proves the lemma.
\end{pf}

Note that the existence of morphism $T(\alpha)$ for any acyclic
fibration $\alpha$ already implies that weakly equivalent algebras
have weakly equivalent tangent Lie algebras. 

\subsection{Standard acyclic cofibrations}

Let $\alpha:A\to B$ be a standard acyclic cofibration. This means
that $\alpha$ is isomorphic to a canonical injection 
$A\to A\coprod F(M)$ where $M\in C(k)$ is a contractible complex
of free $k$-modules and $F(M)$ is the corresponding free algebra. Put

$$ T_{\alpha}=\{\delta\in T_B|\delta(A)\subseteq A\}.$$

Denote by $\kappa_{\alpha}:T_{\alpha}\to T_B$ the natural inclusion.
Note that $T_{\alpha}$ is a dg Lie subalgebra and a Lie algebra map
$\rho_{\alpha}:T_{\alpha}\to T_A$ is defined.

\subsubsection{}
\begin{prop}{T(ac)}The map $\rho_{\alpha}$ is a surjective 
quasi-isomorphism while $\kappa_{\alpha}$ is an injective 
quasi-isomorphism.
\end{prop}
\begin{pf}
Step 1.. Prove that $\rho_{\alpha}$ is surjective. 
Since $\alpha$ is a standard acyclic cofibration, any derivation of $A$
can be trivially extended by zero to a derivation of $B$. 
It will automatically belong to $T_{\alpha}$.

Step 2. One has
$$\Ker\rho_{\alpha}=\{\delta\in T_B|\delta(A)=0\}.$$
Put $B=A\coprod F(M)$ as above. Then
$\Ker\rho_{\alpha}=\CHom(M,B)$ is contractible.

Step 3. Exactly as in Step 3 of~\Prop{T(af)} the diagram
\begin{center}
   \begin{picture}(4,3)
	\put(0,0){\makebox(1,1){$T_{\alpha}$}}
	\put(0,2){\makebox(1,1){$T_A$}}
	\put(2,0){\makebox(2,1){$T_B$}}
	\put(2,2){\makebox(2,1){$\Der_{\alpha}(A,B)$}}
	\put(0.5,1){\vector(0,1){1}}
	\put(3,1){\vector(0,1){1}}
	\put(1,0.5){\vector(1,0){1.4}}
	\put(1,2.5){\vector(1,0){0.7}}

	\put(1,2.5){\makebox(0.7,0.5){$\scriptstyle\alpha_*$}}
	\put(1,.5){\makebox(1.4,0.5){$\scriptstyle\kappa_{\alpha}$}}
	\put(0.5,1){\makebox(0.7,1){$\scriptstyle\rho_{\alpha}$}}
	\put(3,1){\makebox(0.7,1){$\scriptstyle\alpha^*$}}
   \end{picture}
\end{center}
is commutative and therefore we get that $\kappa_{\alpha}$ is also 
a weak equivalence.
\end{pf}

\Prop{T(ac)} allows one to define the map $T(\alpha): T_A\to T_B$
in the homotopy category $\hol(k)$ as 
$$T(\alpha)=\kappa_{\alpha}\circ\rho_{\alpha}^{-1}.$$

\subsubsection{}
\begin{lem}{comp-ac} Let $A\overset{\alpha}{\to}B\overset{\beta}{\to}C$
be a pair of standard acyclic cofibrations. Then one has
$$ T(\beta\circ\alpha)=T(\beta)\circ T(\alpha)$$
in $\hol(k)$.
\end{lem}
\begin{pf}
See~\Lem{comp-af}.
\end{pf}

\subsection{... and their comparison}
\subsubsection{}
\begin{prop}{comp-ac-af}
Let $A\overset{\alpha}{\to}B\overset{\sigma}{\to}C$ be a pair of morphisms
so that $\alpha$ is a standard acyclic cofibrations and 
$\sigma,\sigma\circ\alpha$
are acyclic fibrations. Then one has
$$ T(\sigma)\circ T(\alpha)=T(\sigma\circ\alpha)$$
in the homotopy category $\hol(k)$.
\end{prop}
\begin{pf}
Consider the diagram
\begin{center}
   \begin{picture}(9,9)
	\put(0,4){\makebox(1,1){$T_A$}}
	\put(2,4){\makebox(1,1){$T_{\alpha}$}}
	\put(6,4){\makebox(1,1){$T_{\sigma}$}}
	\put(8,4){\makebox(1,1){$T_C$}}
	\put(4,0){\makebox(1,1){$T_{\sigma\circ\alpha}$}}
	\put(4,2){\makebox(1,1){$T$}}
	\put(4,6){\makebox(1,1){$T_B$}}
	
	\put(2,4.5){\vector(-1,0){1}}
	\put(7,4.5){\vector(1,0){1}}
	\multiput(4.44,2)(0,-0.1){11}{$\cdot$}
	\put(4.5,1){\vector(0,-1){0}}
	\put(3,5){\vector(1,1){1}}
	\put(6,5){\vector(-1,1){1}}
	\put(5,3){\vector(1,1){1}}
	\put(4,3){\vector(-1,1){1}}
	\put(4,1){\vector(-1,1){3}}
	\put(5,1){\vector(1,1){3}}
	
	\put(1,4.5){\makebox(1,0.5){$\scriptstyle\rho_{\alpha}$}}
	\put(3,5.5){\makebox(0.5,0.5){$\scriptstyle\kappa_{\alpha}$}}
	\put(5.5,5.5){\makebox(0.5,0.5){$\scriptstyle\iota_{\sigma}$}}
	\put(7,4.5){\makebox(1,0.5){$\scriptstyle\pi_{\sigma}$}}
	\put(2,2){\makebox(0.5,0.5)
              {$\scriptstyle\iota_{\sigma\circ\alpha}$}}
	\put(7,2){\makebox(0.5,0.5){$\scriptstyle\pi_{\sigma\circ\alpha}$}}
	\put(4.5,1){\makebox(0.5,1.2){$\scriptstyle\mu$}}

   \end{picture}
\end{center}
Here $T=T_{\alpha}\times_{T_B}T_{\sigma}$ identifies with the Lie subalgebra
$$ T=\{\delta\in T_B|\delta(A)\subseteq A\text{ and }\delta(I)\subseteq I\}
$$
where $I=\Ker \sigma$. Then the dotted map $\mu:T\to T_{\sigma\circ\alpha}$
is defined by the formula $\mu(\delta)=\delta|_A$.

Step 1. Let us check that $\mu$ is surjective. Note that $B=A\coprod F(M)$
where $M\in C(k)$ is freely generated over $k$ by a collection of elements
$m_i\in M$ and their differentials $dm_i$. The map $\sigma:B\to C$
is therefore uniquely defined by its restriction $\sigma\circ\alpha$ to $A$
and by its values $c_i=\sigma(m_i)\in C$. Since $\sigma\circ\alpha$ is 
surjective there exist $a_i\in A$ such that $c_i=\sigma\alpha(a_i)$.
Therefore, choosing an appropriate isomorphism between $B$ and 
$A\coprod F(M)$ (the one sending $m_i$ to $m_i-a_i$) we can suppose that
$M$ belongs to $I=\Ker\sigma$.

If now $\delta:A\to A$ is a derivative vanishing on $\Ker\sigma\circ\alpha=
I\cap A$ we can define $\tilde{\delta}\in T$ by the formula
$$ \tilde{\delta}|_A=\delta,\ \tilde{\delta}|_M=0.$$

Step 2. One has 
$$\Ker\mu=\{\delta\in T_B|\delta(A)=0\text{ and }\delta(I)\subseteq I\}
=\CHom(M,I).$$
This is obviously contractible.

Step 3. Now we see that all the arrows in the above diagram are 
quasi-isomorphisms.
Since the diagram commutes, this proves the claim. 
\end{pf}

\subsection{Final steps}

\subsubsection{}
\begin{lem}{inv}
Let $i:A\to B$ be a standard acyclic cofibration and $p:B\to A$ be left 
inverse to $i$
(so that it is acyclic fibration). Then $T(i)=T(p)^{-1}$.
\end{lem}
\subsubsection{}
\begin{lem}{hom.inv}
Let $\alpha,\beta:A\to B$ be two homotopy equivalent acyclic fibrations. Then 
$T(\alpha)=T(\beta)$ in $\hol(k)$.
\end{lem}
\begin{pf}Let 
$$B\overset{i}{\to} B^I\overset{p_0,p_1}{\rra} B$$ 
be a path object and $h:A\to B^I$ be a homotopy connecting $\alpha=p_0\circ f$
 with 
$\beta=p_1\circ f$. Since $A$ is cofibrant we can suppose that $i$ is a 
standard acyclic cofibration. The maps $p_0$ and $p_1$ are both left inverse 
to $i:B\to B^I$, therefore $T(p_0)=T(p_1)$. Present the map $f:A\to B^I$ as
a composition $f=q\circ j$ where $q$ is an acyclic fibration and $j$ is a
standard acyclic cofibration. According to~\Lem{comp-af} 
$T(p_0\circ q)=T(p_1\circ q)$ and then~\Prop{comp-ac-af} ensures that
$T(\alpha)=T(\beta)$.
\end{pf}

Now the main result of this Section follows.

\subsubsection{}
\begin{thm}{T-is-a-functor}
Let $\CO$ be a $\Sigma$-split operad over $k$.
There exists a functor $T:\hoa(\CO)^{\is}\to\hol(k)^{\is}$ from the homotopy
category of $\CO$-algebras and isomorphisms to the homotopy category of dg Lie
$k$-algebras and isomorphisms which assigns to each cofibrant $\CO$-algebra
the dg Lie algebra  $T_A=\Der_{\CO}(A,A)$.
\end{thm}
\begin{pf}
Any quasi-isomorphism $\alpha:A\to B$ can be presented as a composition
$\alpha=p\circ i$ where $p$ is an acyclic fibration and $i$ is a standard 
acyclic cofibration. In this case we set
$$ T(\alpha)=T(p)\circ T(i).$$
To prove the theorem we have to check that if $p\circ i$ and $q\circ j$ are
homotopic with $p,q$ acyclic fibrations and $i,j$ standard acyclic
cofibrations then 
$T(p)\circ T(i)=T(q)\circ T(j)$ --- see the Picture below.
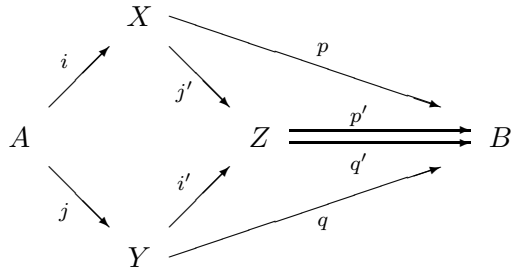
\begin{center}
  \begin{picture}(9,5)
     \put(0,2){\makebox(1,1){$A$}}
     \put(2,0){\makebox(1,1){$Y$}}
     \put(2,4){\makebox(1,1){$X$}}
     \put(4,2){\makebox(1,1){$Z$}}
     \put(8,2){\makebox(1,1){$B$}}

     \put(1,2){\vector(1,-1){1}}
     \put(1,3){\vector(1,1){1}}
     \put(3,1){\vector(1,1){1}}
     \put(3,4){\vector(1,-1){1}}
     \put(5,2.4){\vector(1,0){3}}
     \put(5,2.6){\vector(1,0){3}}
     \put(3,0.5){\vector(3,1){4.5}}
     \put(3,4.5){\vector(3,-1){4.5}}

     \put(1,3.5){\makebox(0.5,0.5){$\scriptstyle i$}}
     \put(1,1){\makebox(0.5,0.5){$\scriptstyle j$}}
     \put(3,3){\makebox(0.5,0.5){$\scriptstyle j'$}}
     \put(3,1.5){\makebox(0.5,0.5){$\scriptstyle i'$}}
     \put(5.3,0.8){\makebox(0.5,0.5){$\scriptstyle q$}}
     \put(5.3,3.7){\makebox(0.5,0.5){$\scriptstyle p$}}
     \put(5.9,1.8){\makebox(0.5,0.5){$\scriptstyle q'$}}
     \put(5.9,2.6){\makebox(0.5,0.5){$\scriptstyle p'$}}
  \end{picture}
\end{center}
Put $Z=X\coprod^AY$. Then the map $A\to Z$ is also a standard acyclic 
cofibration.
Choose $p', q':Z\to B$ so that $p=p'\circ j',\ q=q'\circ i'$. Then obviously 
$p'$ and $q'$
are homotopic. \Prop{comp-ac-af} then says that $T(p)=T(p')\circ T(j')$ and
$T(q)=T(q')\circ T(i')$ which implies 
$$ T(p)\circ T(i)= T(p')\circ T(j')\circ T(i)=T(p')\circ T(i')\circ T(j)=
T(q')\circ T(i')\circ T(j)=T(q)\circ T(j).$$

Theorem is proven.
\end{pf}

\end{document}